\magnification=\magstep1
\noindent

\input amstex
\UseAMSsymbols
\input pictex 
\vsize=23truecm

\NoBlackBoxes
\parindent=18pt
  
   \font\rmk=cmr8      
   
\font\gross=cmbx10 scaled\magstep1


\def\mod{\operatorname{mod}}

\def\Hom{\operatorname{Hom}}
\def\End{\operatorname{End}}
\def\Ext{\operatorname{Ext}}

\def\rad{\operatorname{rad}}

\def\soc{\operatorname{soc}}

\def\D{\operatorname{D}}
\def\endotop{\operatorname{et}}
\def\et{\operatorname{et}}

\def\s{\hfill \square}

\def\und{{\beginpicture
    \setcoordinatesystem units <.1cm,.15cm>
    \multiput{} at -.4 0  2.4 2 /
    \plot 0 2  1 2  1 -.5  2 -.5 /
    \endpicture}}  

\noindent June 2026
	\bigskip
\centerline{\gross Brick chain filtrations. A report.}
	\medskip
\centerline{Claus Michael Ringel}
	\bigskip
{\narrower {\bf Abstract.} We consider the category of finitely generated modules over an
artin algebra $A$. Recall that an object in an abelian category is said 
to be a {\it brick} provided its 
endomorphism ring is a division ring. Simple modules are, of course, bricks, but
in case $A$ is connected and not local, there always do exist bricks which are not simple. 
The aim
of this survey is to focus the attention to filtrations of modules where all factors
are bricks, with bricks being ordered in some definite way, namely according to a
so-called brick chain.

In general, a module category will have many cyclic paths. Recently, Demonet has proposed 
to look at brick chains in order to deal with a very interesting
directedness feature of an arbitrary module category.

The following survey is based on investigations by a large group
of mathematicians. We have singled out some important observations and
have reordered them in order to provide a self-contained 
(and elementary) treatment of the role of bricks in module categories.
(Most of the papers we rely on are devoted to what is called $\tau$-tilting theory,
but for the results we are looking at, there is no need 
to deal with $\tau$-tilting, not even to invoke 
the Auslander-Reiten translation $\tau$ itself). 
\par}
	\bigskip
{\bf Outline.} This is a report on a very important development in the
last 15 years: it focuses the attention to the use of bricks in order to
describe the structure of arbitrary modules over artin algebras.
The report is based on the work of a quite large number of mathematicians,
see section 12.
 We have singled out decisive observations and
have reordered them in order to obtain a self-contained 
and elementary (however still incomplete)
treatment of the role of bricks in module categories.

The first three sections describe the main results presented in the survey, 
they deal with brick chain filtrations and their background. 
Theorem 1.2 and its strengthening 3.2 concern the existence of brick chain filtrations
(and 3.2 includes a corresponding finiteness assertion).
The main tool is the study of torsion classes
and their lower neighbors. 
Theorem 2.3 asserts that finitely generated torsion classes are always generated
by finite semibricks. Theorem 2.8 describes the lower neighbors of the torsion
class generated by a module $M$ in terms of the so-called top bricks of $M$.

Given a brick $B$, we denote by $\Cal E(B)$ the class of all
modules which have a filtration with all factors isomorphic to $B$; these modules
will be said to be {\it homogeneous} of brick type $B$. The brick type of a non-zero
homogeneous module is uniquely determined (see 9.7).
The brick chain filtrations studied in this report concern filtrations
of modules with factors in suitable subcategories $\Cal E(B),$ namely using bricks $B$
which occur in a brick chain. 
The existence of brick chain filtrations is derived from a result 
for neighbor torsion classes.  
Neighbor torsion classes $\Cal T' \subset \Cal T$ come with a label: this is a brick $B$
with the following property: any module $M$ in $\Cal T$ has a submodule $M'$ in
$\Cal T'$ such that $M/M'$ belongs to $\Cal E(B)$, see Theorem 2.7.
	\bigskip\medskip
{\bf 1. All modules have brick chain filtrations.}
	\medskip 
{\bf 1.1.} We deal with an artin algebra $A$; the modules to be considered are usually 9.
left $A$-modules of finite length. 
Given a set $\Cal X$ of modules, let $\Cal E(\Cal X)$ be the class of modules
which have a filtration with all factors in $\Cal X$. If $M_1,\dots,M_m$ are modules, let
$\Cal E(M_1,\dots,M_m) = \Cal E(\{M_1,\dots,M_m\})$ (such a convention is
used throughout the paper in similar situations).

We recall that a {\it brick} is a module whose endomorphism ring is a division ring.
If $B$ is a brick, the modules in $\Cal E(B)$ will be said to be {\it homogeneous}
of brick type $B$.
A finite sequence $(B_1,\dots, B_m)$ is called 
a {\it brick chain,} if all $B_i$ are bricks
and $\Hom(B_i,B_j) = 0$ for $i<j$ (in Appendix 11, we also will deal with
infinite sequences). 
A filtration $0 = M_0 \subset M_1 \subset \cdots \subset M_m = M$ will be called
a {\it brick chain filtration,} provided there is a brick chain
$(B_1,\dots, B_m)$ (its {\it type}) such that $M_i/M_{i-1}$ is homogeneous of brick type $B_i$,
for all $1 \le i \le m$ (the number $m$ is called the {\it length} of the filtration). 
		\medskip
{\bf 1.2. Theorem.} {\it Any module has  brick chain filtrations.} 
	\medskip
The result will be strengthened in 3.2. The proof of Theorem 3.2 is given in section 9.
The assertion 1.2 is also an immediate consequence of Theorem 11.2. 

If all composition factors of a module $M$ are isomorphic, then $(0 \subseteq M)$
is a brick chain filtration of $M$, and obviously the only one;
{\it otherwise, $M$ has at least two different
brick chain filtrations,} as we will show in 10.3.
	\medskip
{\bf 1.3. Some examples of brick chain filtrations.} 
	\smallskip
(1) Let $S_1,\dots,S_n$ be the simple $A$-modules. Obviously, 
$(S_1,\dots,S_n)$ (in any order!) 
is a brick chain. Let us now assume that $\Ext^1(S_i,S_j) = 0$ for
all $i > j.$ If $M$ is an $A$-module $M$, let $M_i$ be 
the maximal submodule of $M$ with all 
composition factors of the form $S_1,\dots,S_i$. 
If $M$ is sincere, then $(M_i)_i$ is a brick chain filtration of type $(S_1,\dots,S_n)$ 
(if $M$ is not sincere, $(M_i)_i$ is  
a brick chain filtration with repetitions, as considered in 10.4).

In particular, recall that $A$ is said to be {\it directed,} provided the
simple modules $S_1,\dots,S_n$ can be ordered in such a way that $\Ext^1(S_i,S_j) = 0$
for all $i \ge j.$ {\it For a directed algebra $A,$ any sincere $A$-module $M$ has a
brick chain filtration of type $(S_1,\dots,S_n)$ such that all 
factors of the filtration are semisimple.} 

(2) Let us consider the special case where $M = B$ itself is a brick. In this case, 
there is, of course, the trivial brick chain filtration $(0 \subset B)$, but if $B$ is not
simple, then there are additional brick chain filtrations (see 10.3). 

Let $A$
be the path algebra of the linearly directed quiver of type $\Bbb A_n$ and $B$ the indecomposable
sincere $A$-module. Then any proper filtration of $B$ is a brick chain filtration, thus
$B$ has $2^{n-1}$ brick chain filtrations. In particular, $B$ 
has brick chain filtrations of length $m$, for any $1 \le m \le n.$

For the quiver with vertices $1,2$, an arrow $1\leftarrow 2$
and a loop $\gamma$ at $1$, with relation $\gamma^t$ for some $t\ge 2,$
the projective cover $M$ of the simple module $2$ is a brick. It 
has the brick chain filtration 
$(0 \subset M_1 \subset M)$ of type $(B_1,B_2)$, where $M_1 = \rad M$, and 
$B_1 = 1, B_2 = 2$.
The relative Loewy length of $M_1$ in $\Cal E(B_1)$ (of course, also its absolute
Loewy length) is $t$.


(3) In contrast to many questions in representation theory, 
looking for brick chain filtrations of modules, it does not seem to be 
helpful to consider first indecomposable modules. Namely, 
brick chain
filtrations of modules $M$ and $M'$ usually do not provide  
a brick chain filtration of $M\oplus M',$ see 10.5. 

(4) (Duality)
Let us denote by $\D$ the usual duality functor.
Given a brick chain filtration $0 = M_0 \subset M_1 \subset \cdots 
\subset M_m = M$ of type $(B_1,\dots,B_m),$ then clearly $\D$ yields a 
corresponding brick chain filtration $(N_i)_i$ of $N = \D M,$ namely
$N_i = \D M/\D M_{m-i}$, for $0 \le i \le m.$ The type of the filtration $(N_i)_i$
is $(\D B_m,\dots ,\D B_1).$ 

(5) Our proof of 1.2 will yield quite special brick chain filtrations,
namely ``torsional'' ones, see section 3.
Let us note already here: if a filtration
$(M_i)_i$ of a module $M$ is torsional, then 
the top of any $M_i$ is generated by the top of $M$. 
Thus, even in the case of a directed algebra, the brick chain filtrations which we
will construct are usually different from the obvious filtrations mentioned in (1).
	\bigskip
{\bf 2. Torsion classes, in particular the finitely generated ones.}
	\medskip
The proof of Theorem 3.2 will be based on the use of torsion classes,
they are essential for all considerations. 
Here, wee recall the
definition and some properties of torsion classes.
	\medskip
{\bf 2.1.}
A class $\Cal T$ of modules is said to be a {\it torsion class} provided $\Cal T$
is closed under factor modules and extensions. The set of all torsion classes 
is a complete lattice; the meet of a set of torsion classes is just the
set-theoretical intersection. Given a class $\Cal X$
of modules, we denote by $T(\Cal X)$ the smallest torsion class which contains
$\Cal X$ (thus, the closure of $\Cal X$ under factor modules and extensions, or,
equivalently, the set-theoretical intersection of all torsion classes containing
$\Cal X$).  The Noether theorems show that
$T(\Cal X)$ is the class
of modules which have a filtration whose factors are factor modules of modules in 
$\Cal X$. The 
{\it torsion submodule} of a module $M$ with respect to the torsion class $\Cal T$
is by definition the largest submodule of $M$ which belongs to $\Cal T.$
Given a module class $\Cal Y$, we denote by ${}^\perp \Cal Y$ the class of all
modules $X$ such that $\Hom(X,Y) = 0$ for all modules $Y$ in $\Cal Y$. It is clear that
${}^\perp \Cal Y$ is closed under factor modules and extensions, thus it is a torsion class. 

A torsion class $\Cal T$ is said to be {\it finitely generated} provided there is
a module $M$ with $\Cal T = T(M).$ Of course, any torsion class $\Cal T$ is the
set-theoretical union of the finitely generated torsion classes contained in $\Cal T$.
	
	\medskip
{\bf 2.2.}
Let $(M_i)_i$ be a brick chain filtration of the module $M$, say of type
$(B_1,\dots,B_t).$ Then $M_i$ is the torsion submodule of $M$ with respect to
the torsion class $T(B_1,\dots,B_i)$ and also the torsion submodule of $M$
with respect to the torsion class ${}^\perp\{B_{i+1},\dots,B_t\}$. Thus, we see:
{\it Given a module $M$ with a brick chain filtration $(M_i)_i$,
the submodules $M_i$ are uniquely determined by the type of the filtration.}
	\medskip
{\bf 2.3.} 
Modules $M, M'$ are defined to be {\it $\Hom$-orthogonal} 
provided $\Hom(M,M') = 0 = \Hom(M',M)$.
A {\it semibrick} is a set of pairwise 
$\Hom$-orthogonal bricks. A torsion class which is generated
by a semibrick is said to be {\it widely generated.}
When we deal with sets of (pairwise non-isomorphic) modules, for example when we consider semibricks, these sets will not necessarily be finite (so that we cannot deal or better
do not want to deal with the corresponding direct sum). 

	\medskip
{\bf Theorem.} {\it For any artin algebra $A$, 
the map $\Cal B \mapsto T(\Cal B)$ provides a bijection between finite 
semibricks and the finitely generated torsion classes.}
	\medskip
The surjectivity of the map asserts that
{\it any finitely generated torsion class is widely generated.}
The injectivity assertion can be extended as follows: the map
$\Cal B \mapsto T(\Cal B)$ is a bijection between arbitrary semibricks and the
widely generated torsion classes, see 2.8.

The proof of Theorem 2.3 is given in 5.6 (the surjectivity of the map), 
and in section 8 (the injectivity of the map). In Section 5, we construct
explicitly an inverse of the map $\Cal B \mapsto T(\Cal B)$, for $T(\Cal B)$ being
finitely generated. Let us outline the construction already here.
	\medskip
{\bf Addendum to Theorem.} 
Given a module $M$, 
we define in 5.4 its ``iterated endotop'' $X = \et^\infty M$; 
this is a factor module of $M$. The indecomposable direct summands of $X$ 
are called the {\it top bricks} of $M$ and we denote by $\Cal B(M)$ the 
set of top bricks of $M$. The surjectivity assertion in Theorem 2.3
can be strengthened as follows (see sections 5 and 8):
{\it Given any module $M$, then $T(M) = T(\Cal B)$ for a uniquely determined semibrick
$\Cal B$, namely the finite semibrick $\Cal B = \Cal B(M)$ of the top bricks of $M$.}
In particular:
{\it $T(M)$ is generated by
a finite semibrick $\Cal B$ whose members are factor modules of $M$.}
	\medskip
{\bf 2.4. Remark.} 
The bijection provided by Theorem 2.3 is of great importance, since it allows to
consider the set of finite semibricks as a partially ordered
set, using the natural partial ordering of the set of torsion classes,
given by set-theoretical inclusion. 
This poset structure on the set of finite semibricks (thus also on the set of bricks)
provides the foundation for the notion of a brick chain as used in Theorem 1.2 and in 
section 11.
	\medskip
{\bf 2.5. Brick finite algebras.} An algebra $A$
is said to be {\it brick finite} provided there are only finitely many  
isomorphism classes of bricks. The algebra $A$ is said to be {\it torsion class finite} 
provided there is only a finite number of torsion classes. 
	\medskip
{\bf Proposition.} {\it For any algebra, the number of finite 
semibricks is equal to the number of finitely generated torsion classes.

An algebra is brick finite iff it is torsion class finite, and
in this case any torsion class is finitely generated.}
	\medskip
Actually, also the converse of the last assertion is true: 
{\it If all torsion classes
are finitely generated, then the algebra is brick finite,} see
12.13. 
	\medskip
Proof of Proposition. For the first assertion, see 2.3.
If $A$ is torsion class finite, then $A$ has only finitely many 
finitely generated torsion classes, thus only finitely many semibricks, 
thus only finitely many bricks. 
Conversely, assume that $A$ is brick finite, thus 
$A$ has only finitely many finitely generated torsion classes.
Given a torsion class $\Cal T$, one can 
start to construct an inclusion chain of finitely generated torsion classes
$\Cal T = \Cal T_0 \subset \Cal T_1 \subset \cdots \subset \Cal T_t\subseteq \Cal T$.
This process will stop after finitely many steps, thus $\Cal T$ is finitely
generated. We see in this way that all torsion classes are finitely generated. 
Thus there are only finitely many torsion classes. 
$\s$
	\medskip
{\bf 2.6. Neighbors.} 
The torsion classes $\Cal T' \subset \Cal T$ will be said to be {\it neighbors}
provided there is no torsion class $\Cal N$ with 
$\Cal T' \subset \Cal N \subset \Cal T$.  If 
$\Cal T' \subset \Cal T$ are neighbor torsion classes, 
$\Cal T'$ is called a {\it lower neighbor} of $\Cal T$ and $\Cal T$ is called
an {\it upper neighbor} of $\Cal T'$. 
	\medskip
{\bf 2.7. Theorem.} {\it If $\Cal T' \subset \Cal T$ are neighbor torsion
classes, then there is a unique brick $B$ in $\Cal T$ such that 
$\Cal T' = \Cal T\cap {}^\perp B.$} 
This brick $B$ is called the {\it label} of the inclusion
$\Cal T' \subset \Cal T.$
	\smallskip
{\it If $\Cal T' \subset \Cal T$ are neighbor torsion classes
with label $B$, then we have:
\item{\rm(1)} $\Cal T' = \Cal T\cap {}^\perp B$ and $\Cal T = T(\Cal T',B)$.
\item{\rm(2)} 
Any module $M$ in $\Cal T$ has a submodule $M'$ in $\Cal T'$
such that $M/M'$ belongs to $\Cal E(B)$.\par}
	\smallskip
{\it Conversely, assume that $\Cal T' \subset \Cal T$ are torsion classes
and that $B$ is a brick in $\Cal T$. If one of the conditions {\rm(1)} or {\rm(2)} is satisfied,
then $\Cal T' \subset \Cal T$ are neighbors with label $B$.}
	\medskip
For the proof, see 7.2. For 
a further  discussion of neighbor torsion
classes, see 7.5. 
	\medskip
Next, we consider the lower neighbors of some torsion classes. 
If $\Cal T$ is a torsion class, we say that $\Cal T$ has {\it sufficiently many
lower neighbors} provided any torsion class $\Cal N$ with $\Cal N \subset \Cal T$
is contained in a lower neighbor of $\Cal T.$ Similarly,  
$\Cal T$ has {\it sufficiently many
upper neighbors} provided any torsion class $\Cal N$ with $\Cal T \subset \Cal N$
contains an upper neighbor of $\Cal T.$ 

	\medskip
{\bf 2.8. Theorem.} {\it A torsion class $\Cal T$ is widely generated iff $\Cal T$ has
sufficiently many lower neighbors.
If $\Cal T$ is widely generated, say $\Cal T = T(\Cal B)$ with $\Cal B$ a semibrick, 
then $B \mapsto \Cal T\cap{}^\perp B$
is a bijection between the elements of $\Cal B$ and the lower neighbors of $\Cal T.$}
	\medskip 
Theorem 2.8 implies that {\it the map
$\Cal B \mapsto T(\Cal B)$ is a bijection between the set of semibricks and the set
of widely generated torsion classes.} For the proof of theorem 2.8, see section 8. 
	\medskip

{\bf 2.9.}  
A torsion class $\Cal T$ is said to be {\it completely join irreducible} provided the
join $\Cal T_*$ 
of the torsion classes properly contained in $\Cal T$ is still properly contained
in $\Cal T$ (and thus $\Cal T_*$ is 
a lower neighbor of $\Cal T$). Note that $\Cal T$ is completely join irreducible iff
$\Cal T$ has a unique lower neighbor and has sufficiently many lower neighbors.
	\medskip
{\bf Corollary.} {\it The map $B \mapsto T(B)$ provides a bijection between the
isomorphism classes of the bricks and the completely join irreducible torsion classes.}
	\medskip
Proof. Theorem 2.3 sends a brick to the torsion class $T(B)$. According to 2.8, 
$\Cal T(B)$ has a unique lower neighbor, namely $\Cal T_* = T(B)\cap{}^\perp B$
and any torsion class properly contained in $\Cal T$ is contained in $\Cal T_*$. This
shows that $T(B)$ is completely join irreducible. 

Conversely, assume that $\Cal T$ is a completely join irreducible torsion class. 
Clearly, $\Cal T$ is finitely generated: Let $M$ be any module in
$\Cal T \setminus \Cal T_*,$ where $\Cal T_*$ is the join of the torsion classes
properly contained in $\Cal T,$ then $\Cal T = T(M).$ Let $B_1,\dots,B_t$ be the
top bricks of $M$, thus $\Cal T = T(M) = T(B_1,\dots,B_t)$.
According to 2.8, $\Cal T$ has $t$ lower neighbors. Since $\Cal T$ is completely
join irreducible, we have $t = 1,$ thus $\Cal T$ is generated by a brick. $\s$
	\medskip 
{\bf 2.10. Warnings.} Let $M$ be a module. 
If $T(M)$ is a finitely generated torsion class and $B$
a top brick of $M$, the lower neighbor torsion class $T(M)\cap {}^\perp B$ is not 
necessarily finitely generated! A typical example will be presented in 2.11.

Also, we have seen in 2.8 that $T(M)$ has only finitely many lower neighbors. 
What about upper neighbors? If $T(M)$ is finitely generated and $\Cal T''$ is an upper
neighbor of $T(M)$, then trivially $\Cal T''$ is again finitely generated, namely
equal to $T(M\oplus N)$, where $N$ is any module in $\Cal T''\setminus T(M).$ 
However, whereas a finitely generated torsion class has only finitely many lower
neighbors, it may have infinitely many upper neighbors. For a typical example, we again
refer to 2.11.
	
	\medskip
{\bf 2.11. An example: The Kronecker algebra.} 
For the benefit of the reader, we want to consider one example in detail,
the Kronecker algebra $A$; this is the path algebra of the quiver with two vertices $1,2$
and two arrows $1\leftleftarrows 2$. 
(It is the usual example which everyone interested in torsion classes of artin algebras has in mind). 

If $x$ is a vertex of a quiver, the simple representation corresponding to $x$ will
also be denoted by $x$; and $P(x)$ and $I(x)$ will denote the projective cover or the
injective envelope of $x$, respectively (provided they exist). 

For the Kronecker algebra $A$ over the field $k$,
there is the
well-known trisection of the indecomposable $A$-modules: there are the preprojective
modules $\Cal P$, the regular modules $\Cal R$ and the preinjective modules $\Cal I$.
This trisection
gives rise to two important torsion classes: the class $T(\Cal I)$ of the
direct sums of preinjective modules, 
and the class $T(\Cal R)$ of 
the direct sums of preinjective and regular modules. {\it Both torsion classes
$T(\Cal I)$ and $T(\Cal R)$ (as many others) are {\bf not} finitely generated}. 
The class $T(\Cal I)$ is the union of a properly 
ascending chain of torsion classes, thus it is 
not finitely generated. 
Note that $T(\Cal I)$ has no lower neighbor, but infinitely many upper
neighbors. The 
class $T(\Cal R)$ is widely generated, namely by the (infinite!) 
semibrick of the
simple regular modules. 
Thus, $T(\Cal R)$ is not finitely generated. Also, $T(\Cal R)$  
has no upper neighbor (but infinitely many
lower neighbors).

In the following picture of the lattice of all torsion classes of $\mod A$,
the sublattice of all torsion classes $\Cal T$
with $T(\Cal I) \subseteq \Cal T \subseteq T(\Cal R)$ has been dotted.
(As we will outline below: The dotted part is uncountable,
even if $k$ is a finite field! In contrast, outside of the
dotted part, there are always just countably many torsion classes.)
$$
{\beginpicture
    \setcoordinatesystem units <.5cm,.46cm>
\multiput{} at -2 0  2 9 /
\multiput{$\bullet$} at 2 0  3 1  3 2  3 3.5  3 5.5  3 7  3 8  2 9  0 4.5 / 
\plot 3 2.5  3 1  2 0  0 4.5  2 9  3 8  3 6.5 /
\multiput{$\vdots$} at 3 3  3 6.4 / 
\setquadratic
\plot 3 3.5  1.3 4.5  3 5.5   4.7 4.5   3 3.5  /

\setlinear
\setshadegrid span <.6mm>
\vshade 1.3 4.5 4.5  <z,z,,> 2 3.8 5.2  <z,z,,> 3 3.5 5.5 <z,z,,> 4 3.8 5.2  
   <z,z,,>  4.7 4.5 4.5  /
\put{$T(1)$} at -.95 4.5
\put{$\{0\}$} at 2.7 -.4
\put{$T(2)$} at 4.05 .9
\put{$T(I(1))$} at 4.45 1.95
\put{$T(\Cal I)$} at 4 3.05
\put{$T(\Cal R)$} at 4.1 5.8
\put{$$} at 3.5 1
\put{$$} at 3.5 1
\put{$T(P(2))$} at 4.55 8.1
\put{$\mod A$} at 3.3 9.4
\endpicture}
$$

If $\Cal X$ is a {\bf non-empty} set of pairwise 
non-isomorphic simple regular Kronecker modules (thus $\Cal X$ is a non-empty semibrick),
then $\Cal T = T(\Cal X)$ is a torsion class 
with $T(\Cal I) \subset \Cal T \subseteq T(\Cal R)$.
Taking also $T(\Cal I)$ itself into account, 
the torsion classes $\Cal T$ with $T(\Cal I) \subseteq \Cal T \subseteq T(\Cal R)$
correspond bijectively to all the subsets of $\Bbb P^1(k)$
(by definition, $\Bbb P^1(k)$  is the union of 
the one element set $\{\infty\}$ and 
the set of monic irreducible polynomials with coefficients in $k$).
Of course, the set of subsets of $\Bbb P^1(k)$ is uncountable. 

For any torsion class $\Cal T$ with 
$T(\Cal I) \subseteq \Cal T \subseteq T(\Cal R)$,
the number of neighbors of $\Cal T$ is always equal to $\max(|k|,\aleph_0),$
in particular, infinite. 
If $R$ is a simple regular module, then $T(R)$ is (of course) finitely generated,
however its unique lower neighbor 
$T(R)\cap {}^\perp R$ is the torsion class $T(\Cal I)$ which is not finitely generated. 
On the other hand, $T(R)$ has infinitely many upper neighbors. Namely, if
$R'$ is a simple regular Kronecker module, not isomorphic to $R$, then
$T(R\oplus R')$ is an upper neighbor of $T(R)$, and there are infinitely many
such modules $R'$. 

We see that for $\mod A$, 
all torsion classes {\bf but one} (namely $T(\Cal I)$)
have sufficiently many lower neighbors. Dually, 
all torsion classes {\bf but one} (namely $T(\Cal R)$) have sufficiently many upper neighbors.
Also,
every torsion class has either two or else infinitely many
neighbors (those with two neighbors are the functorially finite ones, as
mentioned in 12.13). 

For an arbitrary artin algebra, 
the possible numbers of neighbors of torsion classes
are not known. Note that for some algebras (for example for connected 
wild hereditary algebras and for tubular algebras) there do exists torsion classes 
without any neighbor.
	\bigskip
{\bf 3. Torsional brick chain filtrations.}
	\medskip 
In order to strengthen Theorem 1.2, we need an additional notion.
	\smallskip 
{\bf 3.1.} A submodule $U$ of a module $M$ is said to be {\it torsional} provided
$U$ belongs to $T(M).$ 
A filtration $0 = M_0 \subseteq M_1 \subseteq \cdots \subseteq M_m$ will be said to
be {\it torsional} provided $M_{i-1}$ is a torsional submodule of $M_i$, for all
$1 \le i \le m.$ 

If $(M_i)_i$ is a torsional filtration of $M$, then 
$M_{i-1}$ belongs to $T(M_i)$, for all $1\le i\le t$, thus we have the inclusion chain
$0 = T(M_0) \subseteq T(M_1) \subseteq \cdots \subseteq T(M_m) 
= T(M),$ and therefore all the submodules $M_i$ are torsional
submodules of $M$. If a brick chain filtration $(M_i)_i$, say of type $(B_1,\dots,B_m)$
is torsional, then all the bricks $B_i$ belong to $T(M)$, since $B_i$ is a factor module
of $M_i$ and $M_i$ belongs to $T(M).$
	\medskip
{\bf Warning.} 
{\it A brick chain filtration $(M_i)_i$ of a module $M$, with all $M_i$ being torsional
submodules of $M,$ is not necessarily a torsional filtration!} 
It is easy to exhibit modules $M$ with a filtration $(M_i)_i$ which is not torsional, whereas 
all $M_i$ are torsional submodules of $M$:
Consider for example the 
radical square zero algebra $A$ with two simple modules $1,2$, an arrow $1 \leftarrow 2$ and a
loop at $1$, and take the module $M = I(1)$ of length three. 
Since $T(M) = \mod A,$ all submodules of $M$
are torsional. Let $0 \subset M_1 \subset M_2 \subset M$
be the composition series with $M_2/M_1$ isomorphic to $2$.
The submodule $M_1 = 1$ is not a torsional submodule of $M_2$. 
(It is slightly more difficult to construct such a filtration $(M_i)_i$ which is
in addition a brick chain filtration; an example will be given in 9.6 (4).)
	\medskip
{\bf 3.2. Theorem.} {\it Any module has at least one, but only finitely many 
torsional brick chain filtrations. If $M$ has length $m$, the number of torsional brick chain filtrations of $M$ is bounded by $m!$\,.}
	\medskip
The proof will be given in section 9.
As we will see, the torsional brick chain filtrations of a module $M$
can be constructed easily by induction: 
Let $B$ be a top brick of $M$.
Then 
$M$ has a proper submodule $M'$ which belongs to
$T(M)\cap {}^\perp B$, such that $M/M'$ belongs to $\Cal E(B).$
Since $M'$ is a proper submodule of $M$,
by induction there is a torsional brick chain filtration of $M'$, 
say $0 = M_0 
\subset M_1 \subset \cdots \subset M_{m-1} = M'.$ Let $M_m = M.$ Then 
$(M_i)_{0\le i \le m}$ is a torsional brick chain filtration of $M$.   
	\medskip 
{\bf Questions.} Theorem 3.2 asserts that any module $M$ has only finitely many 
torsional brick chain filtrations. 
Usually, $M$ has plenty additional brick chain filtrations
(see, for example, 10.2). {\it Are there modules with
infinitely many brick chain filtrations?} And: 
{\it Is there a module $M$ of length $m$ with more than $m!$ brick chain filtrations?}
	\medskip
{\bf 3.3.}
If $(M_i)_i$ is a torsional brick chain filtration
of type $(B_1,\dots,B_m),$ then by definition all the bricks $B_i$ belong to $T(M)$. 
The last brick $B_m$ is a factor module of $M$, but 
{\it the remaining bricks 
$B_i$ do not have to be factor modules of $M$.} Here is a typical example:
Let $M$ be serial with composition factors going up: $1,2,2,1,2,$ with
torsional brick chain filtration $0 \subset M_1 \subset M$, where $M_1$ is of
length three; here, $M_1$ is not generated by $M$.
	\medskip
{\bf 3.4.} Given a direct sum $X = M\oplus N$, any brick chain filtration $(X_i)_i$ of $X$ 
gives rise to a filtration of $M$ which is, after deleting repetitions,
a brick chain filtration of $M$, we call it the {\it induced} one (see 10.4 (a))
We will show: {\it Any brick chain filtration is
induced from a torsional brick chain filtration,} see 10.4 (b).
	\bigskip

{\bf 4. Some preliminaries.}
	\medskip 
{\bf 4.1. Lemma.} {\it Let $M'$ be a non-zero module in $T(M)$. 
Then $\Hom(M,M') \neq 0.$}
	\medskip
Proof: $M'$ has a filtration $0 = M'_0 \subset M'_1 \subset \cdots \subset M_m = M$, 
where all the factors $M_i/M_{i-1}$ are non-zero factor modules of $M$.
Since $M'_1$ it is a factor module of $M$, we get
a non-zero homomorphism $M \to M'_1 \to M'.$ $\s$
	\medskip
{\bf 4.2. Examples} of non-isomorphic bricks $B', B$ with $B'\in T(B).$ 
According to Lemma 4.1, $\Hom(B,B') \neq 0$. 
(On the other hand, we will see in 6.3 that $\Hom(B',B) = 0.$)
We sometimes will specify modules by a display of the composition factors.
For the following examples, we may deal with a quiver 
with two vertices, labeled $1$ and $2$, with an arrow $1 \leftarrow 2$, and a loop at $1$. 
The display  $\smallmatrix 2 \cr 1 \endsmallmatrix$ stands for a serial module of
length two with socle $1$ and top $2$, and so on \dots.
	\medskip
(1) Let $B = 
\smallmatrix 2 \cr 1 \endsmallmatrix$ and $B' = 2.$  
There is an epimorphism $B \to B'$. (Or, if we want that $B, B'$ have the same support:
Let $B = \smallmatrix 2 \cr 1 \cr 1 \endsmallmatrix$, 
and $B' = \smallmatrix 2 \cr 1 \endsmallmatrix$.)
	\smallskip
(2) Let  $B = \smallmatrix 2 \cr 1  \endsmallmatrix$, 
and $B' = \smallmatrix 2\cr 2 \cr 1 \endsmallmatrix$.
There is a monomorphism $B \to B'$. 
	\smallskip
(3) Let 
$B = \smallmatrix 2 \cr 1 \cr 1 \endsmallmatrix$, 
and $B' = \smallmatrix 2 \cr 2 \cr 1 \endsmallmatrix$. 
There is a non-zero map $B \to B'$, neither epi nor mono.
	\smallskip
The bricks mentioned here are part of a brick chain of the form $\bigl(\, 2,\  
\smallmatrix 2 \cr 2 \cr 1 \endsmallmatrix,\ 
\smallmatrix 2 \cr 1 \endsmallmatrix,\ 
\smallmatrix 2 \cr 1 \cr 1 \endsmallmatrix,\ 1 \, \bigr).$

	\medskip
{\bf 4.3. Lemma.} {\it A non-zero module is a brick iff it has no non-zero proper torsional submodules.}
	\medskip
Proof. Let $M$ be a module. If $M$ is not a brick, there is an endomorphism $f$
of $M$ such that $f(M)$ is a non-zero proper submodule. Since $f(M)$ belongs to
$T(M)$, we see that $f(U)$ is a torsional submodule of $M$.

Conversely, let $U$ be a non-zero proper submodule
which is torsional. Since $U$ belongs to $T(M),$ there is a non-zero submodule $U'$
of $U$ which is a factor module of $M$. We get a non-zero and not invertible
endomorphism $M \to U' \subseteq U \subset M$, thus $M$ is not a brick.
$\s$
	\bigskip
{\bf 5. The endotop and the iterated endotop of a module.}
	\medskip
We are going to show the surjectivity assertion of Theorem 2.3. We need the notion
of the endotop \ $\et M$ \  of a module $M$.
	\medskip
{\bf 5.1. Endotop.} Denote by $E = \End(M)$ the endomorphism
ring of $M$ (operating on the left of $M$), and $\rad E$ its radical. Then 
$(\rad E) M$ is a submodule of $M$ and we define $\et M = M/(\rad E)M$, 
and call it the {\it endotop}
of $M$; by definition, the endotop of $M$ is a factor module of $M$. 

One may define the endotop $\et M$ also as follows: $\et M = M/M'$,
where $M'$ is the sum of the images of the nilpotent endomorphisms of $M$
(in this way, one avoids the question whether one has to look at the ring $\End(M)$
or its opposite, as well as the related one, 
whether functions are written left or write of the argument).
	\medskip
{\bf 5.2. Examples.} (1)
{\it If $M$ is an indecomposable module, $\et M$ may be decomposable.}
For example, let $A$ be a local algebra with radical-square-zero and ${}_AA$ of
length $t\ge 2$. If $M$ is the indecomposable injective module, 
then $\endotop M$ is the direct sum of $t\!-\!1$ copies of the simple module.

(2)
Let $A$ be given by the quiver $Q$ with one vertex and two loops and with relations all paths of length 3 (thus $A$ is a local algebra of dimension 7). There is a 
serial module $M$ of length 3 with $\rad M$ not isomorphic to $M/\soc M.$ 
Then $\et M = M/\soc M$, thus 
{\it $\et M$ is indecomposable of length two, and not a brick,} in particular,
$\et(\et M)$ is a proper factor module of $\et M.$
This leads us below to consider not only $\et$, but the iterations $\et^i,$ see 5.4.
(Instead of $A$, we may consider a proper factor algebra $A'$ of $A$, namely 
the subring $A' = k+J$ of the ring of all $3\times3$-matrices with
coefficients in $k$, where $J$ is the set of
nilpotent upper triangular matrices; let $M = k^3$ be the $A'$-module of 
column vectors.)

(3)	
If $A$ is the Kronecker algebra, and 
$M$ is a regular Kronecker module, then $\et M$ is just the regular top of $M$. 
	\medskip
{\bf 5.3. Proposition.} {\it Let $M$ be a module. 
Then $M$ belongs to $T(\endotop M),$ 
therefore $T(M) = T(\endotop M)$.
The kernel of the canonical map $M \to \et M$ is torsional.}
	\medskip
Proof.
Let $E = \End M$ and 
let $f_1,\dots,f_t$ be a basis of $\rad E$. Let $(\rad E)^m = 0.$
The image of 
$g = (f_i)\:\bigoplus_i M \to M$ is $(\rad E)M = \rad_EM = M_1$ and $\endotop M = M/M_1.$
Let  $M_{j+1} = g(M_j)$ for all $j\ge 0$ with $M_0 = M$. Then $M_m = 0.$  
By induction, all modules $M_j/M_{j+1}$ are generated by $\endotop M.$ This shows that
$T(M) \subseteq T(\et M).$ 
On the other hand, we also have $T(\et M) \subseteq T(M),$ since
$\endotop M$ is a factor module of $M$. 
Thus $M$ and $\endotop M$ generate the same torsion-class.
The kernel $M'$ of the canonical map $M \to \et M$ is by definition the image of the
map $g,$ thus generated by $M$, thus $M'$ belongs to
$T(M).$ $\s$
	\medskip
{\bf 5.4.} We iterate the construction $\et$ and get epimorphisms
$$
 M \to \et M \to \et^2 M \to \cdots.
$$ 
Since $M$ is of finite length, the sequence stabilizes eventually. there is
a non-negative integer $a$ with $\et^a M = \et^{a+1}M$.  In this way,
we get the {\it iterated endotop} $\et^\infty M = \et^a M$ 
(and we have $\et(\et^\infty M) = \et^\infty M$).

	\medskip
{\bf Example.} Let $A$ be a suitable artin algebra with two simple modules $1$ and $2$.
For $n\ge 0,$ let $M[n]$ be a serial module of length $n+2$, 
with composition factors going up: $(1,\dots,1,2,1)$
(starting with $n$ factors of the form $1$).
Then, for $0\le i \le n,$ we have $\et^{i}M[n] = M[n-i]$. For $0\le i < n$, the
module $M[i]$ is not a brick, but $\et^n M[n] = M[0]$ is a brick of the form 
$\smallmatrix 1\cr 2 \endsmallmatrix$.
	\medskip
{\bf 5.5. Proposition.} {\it Let $M$ be a module. The iterated endotop 
$X = \et^\infty M$ is the direct sum of modules which belong to a semibrick $\Cal B(M)$
and $T(M) = T(X) = T(\Cal B(M)).$} The elements of $\Cal B(M)$ are called the
{\it top bricks} of $M$.
{\it  The kernel of the
canonical map $M \to \et^\infty M$ is a torsional submodule of $M$.}
	\medskip
Proof. It is obvious that the iterated endotop of a module is always the direct sum of 
modules which belong to a semibrick,
since the sequence $ M \to \et M \to (\et)^2 M \to \cdots$
stabilizes precisely when $\End(\et^a M)$ is semisimple. Proposition 5.3 yields that
the torsion classes $T(\et^i M)$ are equal, for all $i \ge 0.$ 

The kernel $K$
of the canonical map $M \to \et^\infty M$ has a filtration whose factors are
the kernels $K_i$ of the canonical maps $\et^{i} M \to \et^{i+1} M$, for all $i \ge 0.$
According to 5.3, all modules $K_i$ belong to $T(M)$, thus $K$ belongs to
$T(M).$
$\s$
	\medskip
{\bf 5.6. Corollary.} {\it A torsion class $\Cal T$ is finitely generated iff 
there is a finite semibrick $\Cal B$  with $\Cal T = T(\Cal B).$} 
$\s$
	\medskip
Corollary 5.6 shows that the map $\Cal B \mapsto T(\Cal B)$ from the set of finite
semibricks $\Cal B$ to the set of finitely generated torsion classes is surjective.
This is part of Theorem 2.3.
	\medskip
{\bf 5.7.} {\bf Examples.} 
(1) {\it Let $M$ be an indecomposable module. 
A top brick of $M$ may occur in $\et^\infty M$ with arbitrarily large multiplicity,}
see 5.2 (1). In particular, $\et^\infty M$ may not be indecomposable!

(2) {\it The number of top bricks of an indecomposable module $M$ may be arbitrarily large:}
Consider the $(t\!-\!1)$-subspace quiver, with sink $1$ and sources
$2,3,\dots,t,$ and add a loop at the sink $1$. For the corresponding 
radical-square-zero algebra, the indecomposable injective module $M = I(1)$ satisfies
$\et^\infty M = \et M = 1\oplus 2 \oplus \cdots \oplus t$, thus $M$ has $t$ top bricks.
	\medskip
{\bf 5.8. Warning.} Let $M$ be a module and $(M_i)_i$ a brick chain filtration of
$M$ of brick type $(B_1,\dots,B_m).$ Then, of course, {\it $B_m$ is a factor module of $M$,
but $B_m$ is not necessarily a top brick of $M$.} 
To have an example in mind, just take any
non-simple brick $M$. Then $M$ has just one top brick, namely $M$, whereas given any simple
module $S$ in the top of $M$, there is a brick chain filtration of $M$ of 
type $(B_1,\dots,B_m)$ such that $B_m = S,$ see 10.2.
	\bigskip
{\bf 6. The essential feature: If $B$ is a brick, $({}^\perp B) \und \Cal E(B)$ is a torsion class.}
	\medskip
Given module classes $\Cal X$ and $\Cal Y$, we write $\Cal X\und\Cal Y$ for the
class of all modules $M$ which have a submodule $M'$ in $\Cal X$ such that $M/M'$ belongs
to $\Cal Y$. 

We are going to show: If $B$ is a brick, then
$$
 T({}^\perp B,B)  = ({}^\perp B) \und \Cal E(B). \tag{$*$}
$$
This describes very nicely the torsion class $T({}^\perp B,B)$. 
Actually, there is a corresponding description for some other  
torsion classes $\Cal T \subseteq  T({}^\perp B,B),$ namely the torsion classes
$\Cal T$ with $\Cal T = T(\Cal T\cap{}^\perp B,B)$, see the following Proposition.
Note that the assumption $\Cal T = T(\Cal T\cap{}^\perp B,B)$ means, in particular,
that $B$ belongs to $\Cal T$. 
	\medskip 
{\bf 6.1. Proposition.} {\it If $B$ be a brick and $\Cal T$ a torsion class with
$\Cal T = T(\Cal T\cap{}^\perp B,B)$, then}
$$
 \Cal T = (\Cal T\cap {}^\perp B) \und \Cal E(B).
$$
	\medskip
Let us add: If $0 \to M' \to M \to M/M' \to 0$ is an exact sequence 
with $M'$ in $\Cal T' = \Cal T\cap{}^\perp B$ and $M/M' \in \Cal E(B),$ 
then $M'$ is just the torsion submodule of $M$ with respect to the
torsion class $\Cal T'$, since $\Hom(\Cal T',\Cal E(B)) = 0.$
	\medskip
The description $(*)$ for $\Cal T = T({}^\perp B,B)$, where $B$ is an arbitrary brick,
is a special case of 6.1. 
Namely, for this torsion class $\Cal T,$ we have ${}^\perp B\subseteq \Cal T$, thus 
$\Cal T \cap{}^\perp B = {}^\perp B$, so we have
$T(\Cal T\cap{}^\perp B,B) =
T({}^\perp B,B) = \Cal T$, thus the assumption is satisfied. This means that 
we can conclude that 
$ \Cal T = (\Cal T\cap {}^\perp B) \und \Cal E(B) = ({}^\perp B) \und \Cal E(B),$
as mentioned in $(*)$.
	\medskip
Proof of Proposition. Let $M$ be a module in $\Cal T.$ 
Let $M'$ be a submodule of $M$ which also belongs to $\Cal T$ with $M/M' \in \Cal E(B),$
and minimal with these two properties. We claim that $M'$ belongs to ${}^\perp B,$
thus to $\Cal T\cap {}^\perp B.$ 

Thus, assume for the contrary that there is a non-zero map $f\:M' \to B.$ 
Since $M'$ belongs to $\Cal T,$ there is a filtration 
$0 = M_0 \subset M_1 \subset \cdots \subset M_m = M'$ such that all factors 
$F_i = M_i/M_{i-1}$ are 
factor modules of  $B$ or belong to $\Cal T\cap {}^\perp B$ (in particular, all $F_i$
belong to $\Cal T$). 
Let $s$ be minimal such that $f|M_s$ is non-zero. 
Thus, $f$ vanishes on $M_{s-1}$ and induces a map $\overline f\:M'/M_{s-1}$ with non-zero
restriction to $F_s = M_s/M_{s-1}$. Let us
denote by $u\:F_s \to M'/M_{s-1}$ the inclusion map. Thus, the composition 
$\overline f\cdot u\: F_s \to  B$ is a non-zero map. 

Now $F_s$ is a factor module of some $B$ or belongs to ${}^\perp B.$
Since there is the non-zero map $\overline f\cdot u\: F_s \to  B$,
we see that $F_s$ is a factor module of $B$. 
Also, since $B$ is a brick, 
there is no non-zero map from a proper factor module of $B$ to $B,$ thus we
see that $F_s = B$ and that the composition 
$\overline f\cdot u\:B = M_s/M_{s-1} \subseteq M'/M_{s-1} \to B$
is an isomorphism. This shows that $u$ is a split monomorphism. 
It follows that there is a submodule $M''$ of $M'$ with $M_{s-1} \subseteq M''$,
such that $M_s \cap M'' = M_{s-1}$ and $M_s + M'' = M'$.
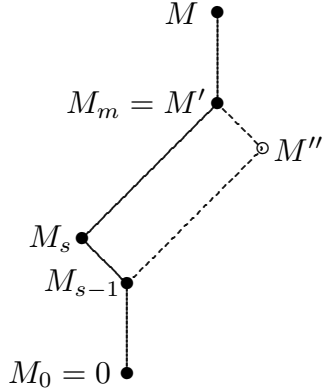
$$
{\beginpicture
    \setcoordinatesystem units <.5cm,.5cm>
\multiput{} at 0 0  4 8 /
\multiput{$\bullet$} at 1 0  1 2  0 3  3 6  3 8 /
\plot 1 0  1 2  0 3  3 6  3 8 /
\put{$\circ$} at 4 5 
\setdashes <.5mm>
\plot 1 2  4 5  3 6 /
\put{$M_s$} at -.7 3 
\put{$M_{s-1}$} at 0 2 
\put{$M_m = M'$} at 1.2 6
\put{$M''$} at 4.8 5
\put{$M_0 = 0$} at -.5 0
\put{$M$} at 2.2 8
\endpicture}
$$
It follows that $M'/M'' \simeq M_s/M_{s-1} = B$, 
 and that 
$M''/M_{s-1}\simeq M'/M_s .$ Since $M/M'$ and $M'/M''$ belong to $\Cal E(B)$, also
$M/M''$ belongs to $\Cal E(B)$. On the other hand, $M''/M_{s-1} \simeq M'/M_s$
has a filtration by factors isomorphic to $F_i$ with $s+1 \le i \le t$
and $M_{s-1}$ has the filtration with factors $F_i$ where $1\le i \le s-1.$
Since all the factors $F_i$ belong to $\Cal T$, also $M''$ belongs to $\Cal T.$

Altogether we see that $M''$ is a submodule of $M$ which belongs to 
$\Cal T$ and such that $M/M' \in \Cal E(B).$ Since $M''$ is a proper submodule
of $M'$, this contradicts the minimality of $M'.$ It follows 
that $M'$ belongs to ${}^\perp B.$ 
 
Since $M'$ belongs to $\Cal T\cap {}^\perp B,$  and $M/M'$ to $\Cal E(B)$,
we see that $M'$ is the torsion submodule of $M$ with respect to the 
torsion class $\Cal T\cap {}^\perp B.$ 
The exact sequence $0 \to M' \to M \to M/M' \to 0$ for an arbitrary module $M$
in $\Cal T$ shows that 
$ \Cal T \subseteq (\Cal T\cap {}^\perp B) \und \Cal E(B)$. On the other hand, we have
$\Cal T\cap {}^\perp B \subseteq \Cal T$, and, since $B\in \Cal T$, also $\Cal E(B)
\subseteq \Cal T:$ This shows the reverse inclusion 
$(\Cal T\cap {}^\perp B)\und \Cal E(B) \subseteq \Cal T$, 
therefore 
$ \Cal T = (\Cal T\cap {}^\perp B) \und \Cal E(B)$. 
$\s$
	\medskip
{\bf Remarks.} If $B$ is a brick which belongs to a torsion class $\Cal T$, 
then it is obvious that 
$T(\Cal T\cap{}^\perp B,B) \subseteq \Cal T,$ but usually, this inclusion will be proper.
see the remark 7.7.

Let us stress that 6.1 can be rephrased as follows: 
The pairs $(\Cal T,B)$, where $B$ is a brick, $\Cal T$ a torsion class, and $\Cal T =
T(\Cal T\cap {}^\perp B,B)$, are just the pairs $(\Cal T,B)$, given by a torsion class
$\Cal T$ with a lower neighbor $\Cal T'$ such that $B$ is the label of the inclusion
$\Cal T' \subset \Cal T$. 
	\medskip
{\bf 6.2. Corollary.} 
{\it Let $B$ be a brick. Let $M$ be a module in $T({}^\perp B,B)$. 
Then any non-zero map $M \to B$ is surjective.}
	\medskip
Proof. Let $M$ be a module in $T(B,{}^\perp B)$ and 
$f\:M \to B$ a non-zero map. The existence of $f$ shows that $M$ does not belong
to ${}^\perp B.$ We have mentioned above that we can use Proposition 6.1 for the torsion
class $\Cal T = T({}^\perp B,B).$ Thus, 
there is a submodule $M'$ of $M$ which belongs to
${}^\perp B$ such that $M/M'$ belongs to $\Cal E(B).$
Since $f$ vanishes on $M'$, we get an induced map
$\overline f\:M/M' \to B$, and $\overline f$ is non-zero. However, any non-zero
map in $\Cal E(B)$ with target $B$ is an epimorphism. Since $\overline f$ is
surjective, also $f$ is surjective.
$\s$
	\medskip	
{\bf 6.3. Corollary.} 
{\it Let $B, B'$ be non-isomorphic bricks, and assume that 
$B'$ is in $T(B)$. Then $\Hom(B',B) = 0$, thus $B' \in T(B)\cap {}^\perp B$.}
	\medskip
Proof. Assume there is a non-zero map $f\:B' \to B$. According to 6.2, the map
$f$ is surjective. Since $B'$ belongs to $T(B)$, we know from 4.1 that
there is a non-zero map $g\:B \to B'$. Since $f$ is surjective, the composition
$gf\:B' \to B \to B'$ is non-zero. Since $B'$ is a brick, this means that
$gf$ is an isomorphism. Thus $f$ is a (split) monomorphism. Altogether we see that
$f$ is bijective, thus $B$ and $B'$ are isomorphic. $\s$
	\bigskip
{\bf 7. Neighbors.}
	\medskip
{\bf 7.1. Lemma.} {\it Let $\Cal T' \subset \Cal T$ be torsion classes. 
Any module $M$ in 
$\Cal T\setminus \Cal T'$ of minimal length is a brick 
and satisfies  $\Cal T' \subseteq {}^\perp M.$}
	\medskip
Proof. Let $M$ be a module in $\Cal T \setminus \Cal T'$ of minimal length.
We form $X = \et^\infty M.$ According to 5.5, we have 
$T(X) = T(M)$, thus also $X$ belongs to $\Cal T \setminus \Cal T'$. There is an
indecomposable direct summand $X'$ of $X$ which belongs to
$\Cal T \setminus \Cal T'$ and, as we know, $X'$ is a brick (one of the top bricks of $M$).
On the other hand, there are epimorphisms $M \to X \to X',$
thus $|X'| \le |X|$. Since we assume that $M$ is of minimal length, we see that
$M = X'$ is a brick. 

In order to see that $\Cal T' \subseteq {}^\perp M,$ consider any homomorphism
$f\: M' \to M,$ with $M' \in \Cal T'.$ Now $f(M')$ belongs to $\Cal T',$ thus
$M/f(M')$ does not belong to $\Cal T'.$ 
Since $M/f(M')$ is a module in $\Cal T\setminus \Cal T'$,
the minimality of $M$ shows that $f(M') = 0.$ 
$\s$
	\medskip
{\bf 7.2. Proof of 2.7.}
First, we assume that $\Cal T' \subset \Cal T$ are neighbors. According to 7.1, there is a brick
$B$ in $\Cal T$ such that $\Cal T' \subseteq {}^\perp B.$ 
Thus, let $B$ be any brick in $\Cal T$ with $\Cal T' \subseteq {}^\perp B.$ 
Then we have
$$
  \Cal T' \subseteq \Cal T\cap {}^\perp B \subset \Cal T
$$
(the proper inclusion is due to the fact that $B$ belongs to $\Cal T$, but not
to ${}^\perp B$).
Since $\Cal T' \subset \Cal T$ are neighbors, we see that 
$\Cal T' = \Cal T\cap {}^\perp B.$ Since $B$ belongs to $\Cal T$ and not to $\Cal T'$,
we also have $\Cal T = T(\Cal T',B).$ This shows that (1) is satisfied.

Next, we show (2): 
any module $M$ in $\Cal T$ has a submodule $M'$ in $\Cal T'$
such that $M/M'$ belongs to $\Cal E(B)$.
It follows from $\Cal T' = \Cal T\cap {}^\perp B$ and $\Cal T = T(\Cal T',B)$ that
$\Cal T = T(\Cal T\cap{}^\perp B,B),$ therefore 
we can use 6.1: For any module $M\in \Cal T,$ there is a submodule $M'$
of $M$ which belongs to $\Cal T\cap {}^\perp B$ with $M/M'$ in $\Cal E(B).$
	\smallskip
Let us show that $B$ is the unique brick in $\Cal T$ with $\Cal T' = \Cal T\cap {}^\perp B.$
Thus, let $C$ be any brick in $\Cal T$ with $\Cal T' = \Cal T\cap {}^\perp C.$
As we have seen, this implies that any module $M$ in $\Cal T$ has a submodule 
$M'$ in $\Cal T'$ with $M/M' \in \Cal E(C).$ 
Now $C$ cannot belong to $\Cal T'$, since otherwise we would have
$\Cal T \subseteq \Cal T'$. The module $C$ has a submodule $C'$ in $\Cal T'$ with $C/C'$
in $\Cal E(B)$. Since $C$ does not belong to $\Cal T'$, we have
$C/C' \neq 0$, therefore $C$ maps onto $B$. 
On the other hand, $B$ has a submodule $B'$ in 
$\Cal T'$ with $B/B'\in \Cal E(C)$. Since $B$ is not in $\Cal T'$, we see that $B/B'$
is non-zero, thus $B$ maps onto $C$. This shows that $C = B.$
	\medskip
Now let $\Cal T' \subset \Cal T$ be any inclusion 
of torsion classes and let $B$ be a brick.
	
First, let us assume that condition (2) is satisfied, thus
any module $M$ in $\Cal T$ has a submodule $M'$ in $\Cal T'$
such that $M/M'$ belongs to $\Cal E(B)$.
Then clearly $\Cal T \subseteq T(\Cal T',B).$ 
In order to show that $\Cal T'\subset \Cal T$
are neighbors, let $\Cal T''$ be a torsion class with 
$\Cal T' \subset \Cal T'' \subseteq \Cal T.$
We claim that $B$ belongs to $\Cal T''.$ 
Let $M$ be a module in $\Cal T'' \setminus \Cal T'$.
According to (2), there is a submodule $M'$ of $M$ in $\Cal T'$ such that
$M/M'$ belongs to $\Cal E(B)$. Since $M$ does not belong to $\Cal T$,
we see that $M/M' \neq 0.$ By definition,
$M/M'$ belongs to $\Cal E(B),$ thus it has a factor module isomorphic to $B$.
Since $M$ belongs to $\Cal T''$, also its factor module
$B$ belongs to $\Cal T''$. Thus we see that $T(\Cal T',B) \subseteq \Cal T''.$
Altogether, we have $T(\Cal T',B) \subseteq \Cal T'' \subseteq \Cal T \subseteq 
T(\Cal T',B)$, therefore $\Cal T'' = \Cal T.$ 

We also claim that $\Cal T' = \Cal T\cap {}^\perp B$. 
The inclusion $\Cal T' \subseteq {}^\perp B$  follows from Lemma 7.1, since
$B$ is a module of minimal length in $\Cal T\setminus \Cal T'.$
(Namely, if $X$ is any module in $\Cal T\setminus \Cal T'$,
then it has a proper submodule $X'\in \Cal T'$ such that $X/X'$ belongs to
$\Cal E(B).$ 
But then $B$ is a factor module of $X/X'$, thus of $X$.)
Since $\Cal T' \subseteq {}^\perp B,$ we have $\Cal T' \subseteq \Cal T\cap {}^\perp B \subset
\Cal T$. But $\Cal T' \subset \Cal T$ are neighbors, thus 
$\Cal T' = \Cal T\cap {}^\perp B$. This shows that the label of 
$\Cal T' \subset \Cal T$ is $B$.

Second, let us assume that (1) is satisfied: 
$\Cal T' = \Cal T\cap {}^\perp B$ and $\Cal T = T(\Cal T',B)$.
Then we have $\Cal T = T(\Cal T\cap{}^\perp B,B),$ thus we can apply 6.1
and see that also condition (2) is satisfied. As we have shown, this implies that 
$\Cal T' \subset \Cal T$
are neighbors with label $B$.
$\s$
	\bigskip
{\bf 7.3. Proposition.} 
{\it Let $B$ be a brick and $\Cal X \subseteq {}^\perp B$. Then
$T(\Cal X,B)\cap {}^\perp B \subset T(\Cal X,B)$, and these are neighbors with label $B$.}
	\medskip
Proof. We write $T(\Cal X,B)_B = T(\Cal X,B)\cap {}^\perp B.$
Now $T(\Cal X,B)_B \subseteq T(\Cal X,B)$, and this inclusion is
proper since $B$ does not belong to ${}^\perp B.$ Assume that there is a torsion
class $\Cal T$ such that 
$T(\Cal X,B)_B \subset \Cal T \subseteq T(\Cal X,B)$.
Since $T(\Cal X,B)_B \subset \Cal T$, there is a module $M\in \Cal T$ which does not
belong to ${}^\perp B.$ Thus, there is a non-zero map $f\:M \to B.$ Since $M$ belongs to
$T({}^\perp B,B)$, we can apply Corollary 6.2. We see that $f$ is surjective, thus $B$
belongs to $\Cal T.$ Of course, also $\Cal X \subseteq \Cal T.$
Therefore $T(\Cal X,B) \subseteq \Cal T.$ This shows that
$T(\Cal X,B) = \Cal T.$ Thus, $T(\Cal X,B)_B \subset T(\Cal X,B)$ are neighbors.
By definition, the label is $B$.
$\s$
	\medskip
{\bf 7.4. Corollary.}
{\it Let $\Cal T' \subset \Cal T$ be torsion
classes. Then there are bricks $B$ in $\Cal T$ such that $\Cal T'\subseteq {}^\perp B.$
If $B$ is a brick with $\Cal T'\subseteq {}^\perp B,$
let $\Cal N = 
T(\Cal T',B)$ and $\Cal N' = \Cal N\cap {}^\perp B,$ then we have
$$
  \Cal T' \subseteq \Cal N' \subset \Cal N \subseteq \Cal T
$$
and the torsion classes $\Cal N' \subset \Cal N$ are neighbors with label $B$.}
	\medskip
Proof: It is trivial that $\Cal T' \subseteq \Cal N' \subseteq \Cal N \subseteq \Cal T.$
We use 7.3 with $\Cal X = \Cal T'$ in order to see that 
$\Cal N' \subset \Cal N$ are neighbor torsion classes with label $B$.
$\s$
	\medskip
{\bf 7.5. Remark.} For a torsion class $\Cal T$, there is the corresponding {\it torsionfree 
class} $\Cal T^\perp$ (by definition, for any class $\Cal X$ of modules, 
$\Cal X^\perp$ is defined as the class of all modules $Y$ 
with $\Hom(X,Y) = 0$ for all $X\in \Cal X$);
the pair $(\Cal T,\Cal T^\perp)$ is called a {\it torsion pair}. 
Subcategory of the form $\Cal E(B)$ with $B$ a brick will be said to be {\it homogeneous
categories} (of brick type $B$).

	\medskip
{\bf Proposition.} {\it Torsion classes $\Cal T' \subseteq \Cal T$ are neighbors iff
$\Cal T\cap (\Cal T')^\perp$ is a homogeneous category. And in this case, the brick type
of $\Cal T\cap (\Cal T')^\perp$ is the label of the inclusion 
$\Cal T' \subset \Cal T$.}
	\medskip
The proposition shows the symmetry between torsion classes and torsionfree classes
when dealing with labels: In particular, the labeling of neighbor torsion classes 
yields a corresponding labeling of neighbor torsionfree classes.
	\medskip
Proof of proposition.
Let $\Cal T' \subseteq \Cal T$ be torsion classes, and let
$\Cal F' = (\Cal T')^\perp.$ First, let us assume that $\Cal T\cap \Cal F'$ is homogenous,
say of type $B$, with $B$ a brick.
Since $B$ belongs to $\Cal F'$, we see that $B$ cannot belong to $\Cal T'$,
thus we have $\Cal T' \subset \Cal T$. 

In order to show that $\Cal T'\subset \Cal T$
are neighbors, 
let $\Cal T''$ be a torsion class with $\Cal T' \subset \Cal T'' \subseteq \Cal T.$
We claim that $B$ belongs to $\Cal T''.$ 
Let $M$ be a module in $\Cal T'' \setminus \Cal T'$ and let $M'$ be its torsion module
with respect to the torsion class $\Cal T'$. Since $M$ does not belong to $\Cal T'$,
we see that $M/M' \neq 0.$ By definition,
$M/M'$ belongs to $\Cal F'.$ Since $M$ belongs to $\Cal T$, also its factor module
$M/M'$ belongs to $\Cal T$, therefore $M/M'$ belongs to $\Cal T \cap \Cal F' =
\Cal E(B)$. As a non-zero module in $\Cal E(B)$, the module $M/M'$ has $B$ as
a factor module, therefore $B$ belongs to $\Cal T''.$ 

Now, let $X$ be any module in $\Cal T$. 
Let $X'$ be the torsion submodule of $X$ with respect to $\Cal T'$, thus
$X/X'$ belongs to $\Cal F'$. Also, $X/X'$ is a factor module of $X\in \Cal T$, thus
$X/X'$ is in $\Cal T$, therefore in $\Cal T\cap \Cal F' = \Cal E(B)$. 
This shows that $X$ has a filtration with modules in $\Cal T'$ and in $\Cal E(B)$,
 thus $X$ belongs to $\Cal T''$.
This shows that $\Cal T'' = \Cal T.$
	\medskip
Conversely, let $\Cal T' \subset \Cal T$ be neighbors, say with label $B$.
Then $B$ belongs to $\Cal T$. Also, $\Cal T' = \Cal T\cap{}^\perp B \subseteq {}^\perp B$
shows that $B \in (\Cal T')^\perp = \Cal F'.$ Since 
$\Cal T \cap \Cal F'$ is closed under extensions, $\Cal E(B) \subseteq
\Cal T \cap \Cal F'$. 
It remains to show that 
$\Cal T \cap \Cal F' \subseteq \Cal E(B).$ 
	
First, we claim that $\Cal T' = \Cal T\cap {}^\perp B$ and that $T(\Cal T',B) = \Cal T$.
Namely, since $\Cal T' \subset \Cal T$ are neighbors and 
$\Cal T' \subseteq \Cal T\cap {}^\perp B \subset \Cal T$, we have
$\Cal T' = \Cal T\cap {}^\perp B$. 
On the other hand, since $\Cal T' \subset \Cal T$ and $B \in \Cal T,$ 
we have $T(\Cal T',B)\subseteq \Cal T$.
Since $B$ does not belong to $\Cal T'$, we have $\Cal T' \subset T(\Cal T',B)$.
Altogether, we have
$\Cal T' \subset T(\Cal T',B) \subseteq \Cal T$. Thus 
$T(\Cal T',B) = \Cal T$, since $\Cal T' \subset \Cal T$ are neighbors.

Since $\Cal T= T(\Cal T\cap{}^\perp,B)$ we can use propositon 6.1. It asserts:
If $M$ belongs to $\Cal T$ and $M'$ is the torsion submodule of $M$
with respect to $\Cal T\cap {}^\perp B$, then $M/M'$ belongs to $\Cal E(B)$. 
Now assume that $M$ belongs to $\Cal T \cap \Cal F'$. Since $M \in \Cal F'$,
the torsion submodule $M'$ of $M$ with respect to $\Cal T\cap {}^\perp B = \Cal T'$
is zero. Thus,  $M = M/M'$ belongs to $\Cal E(B)$. 
$\s$
	\medskip
{\bf 7.6. The brick chains explained in terms of neighbor torsion classes.} 
Let $\Cal T' \subset \Cal T$ be neighbors with label $B$.
Then we have on the one hand: $B$ belongs to $\Cal T$ and not to
$\Cal T'$. On the other hand, for every module $M$ in $\Cal T'$, 
in particular for the bricks in
$\Cal T'$, we have $\Hom(M,B) = 0.$ 
Thus we obtain in this way the $\Hom$-condition which is used in the definition of a
brick-chain: If $\Cal T_1 \subset \Cal T_2 \subseteq \Cal T_3 \subset \Cal T_4$
is a chain of torsion classes with $\Cal T_1 \subset \Cal T_2$ as well as 
$\Cal T_3 \subset \Cal T_4$ being neighbors, and $B$ is the label for 
$\Cal T_1 \subset \Cal T_2$, whereas $B'$ is the label for $\Cal T_3 \subset \Cal T_4$, 
then $\Hom(B,B') = 0.$
	\medskip
{\bf 7.7. Remark.} If $\Cal T' \subset \Cal T$ are neighbors with
label $B$, we have both  
$T(\Cal T\cap{}^\perp B,B) = \Cal T$
and $T(\Cal T'\cap{}^\perp B,B) = \Cal T'.$ 
In general, for arbitrary torsion classes $\Cal T$ and $\Cal T'$, and $B$ a brick,
there are the obvious inclusions 
$T(\Cal T\cap{}^\perp B,B) \subseteq  \Cal T$ 
provided $B$ belongs to $\Cal T$, as well as
$\Cal T' \subseteq T(\Cal T',B)\cap{}^\perp B$ provided $\Cal T' \subseteq {}^\perp B.$
But both inclusions are usually proper inclusions. We are obliged to the referee for pointing
out the following example. 

Let $A$ be the path algebra of the $\Bbb A_2$-quiver $1\leftarrow 2$ and $B$ the indecomposable
module of length two. Here, we have ${}^\perp B = T(2) \subset T(B)$. 
For $\Cal T = \mod A$, we have 
$T(B) = T(\Cal T\cap{}^\perp B,B) \subset  \Cal T = \mod A.$ 
For $\Cal T' = \{0\}$, we have $\{0\} = \Cal T' \subset T(\Cal T',B)\cap{}^\perp B = T(2).$

In particular, given an inclusion $\Cal T' \subset \Cal T$ and a brick $B$,
the conditions $\Cal T' \subseteq {}^\perp B$ and $\Cal T = T(\Cal T',B)$
do not imply that $\Cal T' \subset \Cal T$ are neighbors. 
	\bigskip
{\bf 8. Widely generated torsion classes.}
	\medskip
We are going to prove 2.8. If $\Cal B$ is a semibrick and $B\in
\Cal B,$ we write $T(\Cal B)_B = T(\Cal B)\cap {}^\perp B.$
	\medskip
{\bf 8.1. Lemma.} {\it Let $\Cal B$ be a semibrick, and $\Cal T = T(\Cal B).$ 
If the torsion class $\Cal T'$ is properly contained in $\Cal T$, then there is
$B \in \Cal B$ with $\Cal T' \subseteq \Cal T_B$
and such that $\Cal T_B \subset \Cal T$ are neighbors.}
	\medskip
Proof. Since $\Cal T'$ is properly contained in $\Cal T$, there is a brick $B \in \Cal B$
which is not contained in $\Cal T'$. Let $\Cal B' = \Cal B \setminus \{B\}.$
Since $\Cal B$ is a semibrick, we have $\Cal B' \subseteq {}^\perp B.$
According to 7.3, $\Cal T_B \subset \Cal T$ are neighbors. 

Also, we claim that $\Cal T' \subseteq {}^\perp B.$ Namely, if $f\:M \to B$ is
a non-zero homomorphism with $M \in \Cal T'$, then 6.2 asserts that $f$ is surjective,
thus $B \in \Cal T',$ a contradiction. 
It follows that $\Cal T' \subseteq {}^\perp B$, thus 
$\Cal T' \subseteq \Cal T\cap {}^\perp B = \Cal T_B.$ 
$\s$
	\medskip
{\bf 8.2. Proof of Theorem 2.8.} First, let $\Cal B$ be a semibrick and $\Cal T = 
T(\Cal B).$ According to 8.1, $\Cal T$ has sufficiently many lower neighbors,
namely the torsion classes $\Cal T_B$ with $B\in \Cal B.$
Also, the map $B \mapsto \Cal T_B$ from $\Cal B$ to the set of lower neighbors of $\Cal T$
is surjective. On the other hand, this map is injective by the unicity of the label.

Conversely, let $\Cal T$ be a torsion class with sufficiently many lower neighbors.
Let $\Cal B$ be the set of labels of the lower neighbors. Then $\Cal B$ is a subset
of $\Cal T$, thus $T(\Cal B) \subseteq \Cal T.$
Let us assume that $T(\Cal B) \subset \Cal T.$ Since $\Cal T$ has sufficiently many lower
neighbors, there is a lower neighbor $\Cal T'$ of $\Cal T$ such that
$T(\Cal B) \subseteq \Cal T'.$ Let $B$ be the label of the inclusion
$\Cal T' \subset \Cal T$. Then $B\in \Cal B$. Now $\Cal T' = \Cal T\cap {}^\perp B
\subseteq {}^\perp B.$ Thus we have $B \in T(\Cal B) \subseteq \Cal T' \subseteq 
{}^\perp B$, a contradiction.
$\s$
 	\medskip
{\bf 8.3. Corollary.} {\it If $\Cal B$ is a semibrick.
Then $T(\Cal B)$ is finitely generated iff $\Cal B$ is finite.}
	\medskip
Proof. If $\Cal B$ is finite, then, of course, $T(\Cal B)$
is finitely generated. Conversely, assume that $T(\Cal B)$
is finitely generated. By definition, there is a module $M$ with 
$T(\Cal B) = T(M).$
According to Theorem 2.3, there is a finite semibrick $\Cal B'$
with $T(M) = T(\Cal B').$ According to 2.6 and 2.8, we have $\Cal B =
\Cal B'$, thus $\Cal B$ is finite. $\s$

	\bigskip
{\bf 9. Torsional brick chain filtrations.}
	\medskip
We are going to prove Theorem 3.2.
	\medskip 
{\bf 9.1. Proposition.}
{\it Let $B$ be a top brick of the module $M$. Then 
$M$ has a proper submodule $M'$ which belongs to
$T(M)\cap {}^\perp B$, such that $M/M'$ belongs to $\Cal E(B).$}
	\medskip

Proof. We want to use 6.1 for $B$ and the torsion class $T(M)$. We have to show
that $T(M) = T(T(M)\cap{}^\perp B,B)$. Of course, 
$T(T(M)\cap{}^\perp B,B) \subseteq T(M),$ since $B$ is a factor module of $M$.
For the reverse inclusion, let $\Cal B'$ be the set of 
top bricks different from $B$. According to 5.5, we have $\Cal B' \subseteq {}^\perp B$
and $T(M) = T(\Cal B',B).$ As a consequence, we have 
$\Cal B' \subseteq T(M)\cap {}^\perp B.$ Therefore
$T(M) = T(\Cal B',B) \subseteq 
T(T(M)\cap{}^\perp B,B)$.
$\s$
	\medskip
{\bf 9.2. Corollary.} {\it Let $B$ be a top brick of the module $M$. Then 
there is a torsional brick chain filtration $(M_i)_i$ of $M$
of some brick type $(B_1,\dots,B_m)$ with $B_m = B.$}
	\medskip
Proof by induction. According to 9.1, there is a proper submodule $M'$ of $M$
which belongs to
$T(M)\cap {}^\perp B$, such that $M/M'$ belongs to $\Cal E(B).$
If $M' = 0,$ then $(0 \subset M)$ is a brick chain filtration of $M$ of
type $(B)$, and this filtration is of course torsional.
Otherwise, by induction, $M'$ has a torsional brick chain filtration $(M_i)_i$
say of type $(B_1,\dots,B_{m-1})$. We put $M_{m-1} = M'$ and $M_m = M.$
Then the filtration $(M_i)_{0\le i \le m}$ is the required torsional brick chain
filtration.
$\s$
	\medskip
Thus any module has at least one torsional brick chain filtration.
	\medskip
{\bf 9.3. Lemma.} {\it Let $M$ be a module, $B$ a brick. Assume that
$M$ has a proper torsional submodule $Y$ in $T(M)\cap {}^\perp B$ such that
$M/Y$ belongs to $\Cal E(B)$. Then $B$ is a top brick of $M$
(and $Y$ is the torsion submodule of $M$ with respect to the torsion class
$T(M)\cap{}^\perp B$).} 

	\medskip
Proof. First, we show that $T(M) = T(Y,B)$. Since $Y$ is a proper submodule
of $M$, we see that $M/Y$ is a non-zero module in $\Cal E(B)$, thus it has
a factor module isomorphic to $B$. Since $B$ is a factor module of $M$, we know that $B$
belongs to $T(M).$ Also, by assumption, $Y$ belongs to $T(M)$. Thus
$T(Y,B) \subseteq T(M).$ On the other hand, $M$ has a filtration with
factors of the from $Y$ and $B$, thus $T(M) \subseteq T(Y,B).$

Next, we calculate the iterated endotop of $Y\oplus B.$ 
We calculate inductively $\et^a (Y\oplus B)$ for all $a \ge 0.$
We claim that $\et^a (Y\oplus B) = Y_a\oplus B$, where $Y_a$ is a factor
module of $Y$ with $\Hom(Y_a,B) = 0.$ 
For $a = 0,$ we put $Y_a = Y.$ Assume we have 
$\et^a (Y\oplus B) = Y_a\oplus B$, where $Y_a$ is a factor
module of $Y$ with $\Hom(Y_a,B) = 0.$ 
Since $\Hom(Y_a,B) = 0,$ the radical maps in the endomorphism ring of
$Y_a\oplus B$ map into $Y_a$. If $U_a$ is the sum of these images, then
$\et^a (Y\oplus B) = Y_{a+1}\oplus B$ with $Y_{a+1} = Y_a/U_a.$ 
Also, we have $\Hom(Y_{a+1},B) = 0$, since any non-zero homomorphism $Y_{a+1} \to B$
would yield a non-zero homomorphism $Y_a \to Y_{a+1} \to B$.
Since we deal with modules of finite length, there is some $a$ such that $U_a = 0,$
and therefore $\et^\infty(Y\oplus B) = Y_a\oplus B$.
This shows that $B$ is a top brick of $Y\oplus B$.

Since $T(M) = T(Y,B)$, we know from 2.8 that the top bricks of $M$ are just the
top bricks of $Y\oplus B.$ Thus $B$ is a top brick of $M.$
$\s$
	\medskip
{\bf 9.4. Corollary.} {\it Let 
$(M_i)_i$ be a torsional brick chain filtration of $M$
of brick type $(B_1,\dots,B_m)$. Then $B_m$ is a top brick of $M$
and $M_{m-1}$ is the torsion submodule of $M$ for the torsion class 
$T(M)\cap {}^\perp B_m$.}
	\medskip
Proof. We apply Lemma 9.3 to $Y = M_{m-1}$ and $B = B_m$. $\s$
	\medskip
{\bf 9.5. Finiteness of the number of torsional brick chain filtrations.} 
For any module $M$, let $\phi(M)$ be the number
of torsional brick chain filtrations of $M$. We show by induction that $\phi(M) \le m!$\,,
where $m$ is the length of $M$. Of course, $\phi(0) = 1.$ 
Now, let $M$ be a non-zero module. There are at most $m$ top bricks $B^{(1)},\dots,
B^{(t)}$. For any top brick $B^{(i)}$, let $M^{(i)}$ be the torsion submodule of $M$ with
respect to the torsion class $T(M)\cap {}^\perp B^{(i)}$. Then $M^{(i)}$ has length at most
$m-1$, thus, by induction, $\phi(M^{(i)}) \le (m-1)!\,.$ Therefore 
$
 \phi(M) = \sum_i \ \phi(M^{(i)}) \le m\cdot (m-1)! = m!\,.
$
Namely, according to 9.4 any torsional brick chain filtrations of $M$ 
is obtained from a torsional brick chain filtration of $M^{(i)}$ by adding the
inclusion $M^{(i)} \subset M.$
$\s$
	\medskip
{\bf Remark.} The basic semisimple modules show that the bound $\phi(M) \le m!$ is optimal. 
	\medskip
{\bf 9.6. Some examples.} 
	\smallskip
(1) If $M$ is homogeneous of brick type $B$, then the only torsional brick chain filtration
of $M$ is $(0 \subseteq M).$) 
If the brick $B$ is not simple, and $M$ is non-zero, then $M$ has 
brick chain filtrations which are not torsional, see 10.3. 

(2) Let $A$ be a Nakayama algebra and $M$ an indecomposable
module. Then {\it $M$ has precisely one torsional brick chain filtration,
and this filtration has length at most two.} Proof. 
Let $B$ be the (unique) top brick of $M$. Let $U$ be the
smallest submodule of $M$ such that $M/U$ belongs to $\Cal E(B)$.
If $U = 0$, then $(0 \subset M)$ is the unique torsional brick chain filtration.
If $U \neq 0$, then $(0 \subset U \subset M)$ is the unique torsional 
brick chain filtration.
$\s$

It follows: {\it If the support of $M$ has cardinality $n$, then 
$M$ has at least $2^{n-1}$ brick chain filtrations.} 
We may assume that $M$ is faithful. If $A$ is a directed algebra, then see 1.6 (2).
We 
assume that the quiver of $A$ has arrows $(i\!-\!1)  \leftarrow i $ for $2 \le i < n$
and $n \leftarrow 1$, and that the top of $M$ is $n$. Let $1 \le b_1 < b_2 < \cdots < b_t =n$
be a sequence of integers (the number of such sequences $(b_i)_i$ is $2^{n-1}$).
For $1 \le i \le t$, 
let $M_i$ be the maximal submodule of $M$ with top $b_i$. 
Thus $0 \subset M_1 \subset \cdots \subset M_t = M$ is a filtration of $M$, and we write 
$B_i = M_i/M_{i-1}$ for $2\le i \le t$. We now
use the torsional brick chain filtration of $M_1$:
Let $B_1$ be the top brick of $M_1$. 
Let $U$ be the smallest submodule of $M_1$ such that 
$M_1/U$ belongs to $\Cal E(B_1)$.
If $U = 0$, then 
$0 \subset M_1 \subset \cdots \subset M_t$ is a brick chain filtration of $M$ of type
$(B_1,B_2,\dots,B_t).$ 
If $U \neq 0$, then 
$0 \subset U \subset M_1 \subset \cdots \subset M_t$ is a brick chain filtration of $M$
of type $(U,B_1,\dots,B_t)$.
$\s$

(The bound $2^{n-1}$ is not optimal: If $A$ is a cyclic Nakayama algebra with $n = 3$, any
indecomoposable module of length 4 has five brick chain filtrations.) 

(3) Duality. 
We obtain from a brick chain filtration $(M_i)_i$
of $M$ a brick chain filtration for $\D M,$ see 1.3 (4). The 
dual filtration of a torsional brick chain filtration may not be torsional.
As an example, let $A$ be a connected Nakayama 
algebra with $n = 2$ and let $M$ be an indecomposable module of 
length three, let $U$ be its socle. 
Then $M$ has the brick chain
filtration $(0 \subset U \subset M)$. This filtration is torsional, but
the dual filtration is not. 
	
There are brick chain filtrations $(M_i)_i$ such that  
neither $(M_i)_i$ nor the dual filtration is torsional. 
As an example, take a connected Nakayama algebra with $n = 3$ and let
$M$ be indecomposable of length $4$. Let $U$ be the submodule of $M$ of length $2$.
Then $(0 \subset U \subset M)$ is a brick chain filtration, neither
this filtration nor its dual is torsional.

(4) We have mentioned in 3.1 that 
a brick chain filtration $(M_i)_i$ of a module $M$, with all $M_i$ being torsional
submodules of $M,$ may not be a torsional filtration. Here is an example:
Let $Q$ be the quiver with vertices $1,2,3$ and with arrows
$1\to 2,\ 2 \to 1,\ 2 \to 3$, such that the path
$1 \to 2 \to 1 \to 2 \to 1$ is a zero relation. 
The projective module $P(2)$ has a submodule $U$ isomorphic to $3$ such that
$M = (\rad P(2))/U$ is the direct sum of a serial module $V$ with factors $3,2,1,2,1$
(going up) and a copy of $3$. Take as  $M_1$ and $M_2$ the submodules of $V$ of
length two and four, respectively. Then $0 = M_0 \subset M_1 \subset M_2 \subset  M$
is a brick chain filtration of type $(M_1,M_2/M_1,M/M_2)$. Both bricks $M_1$ and $M_2/M_1$ 
are factor modules of  $M$. It follows that 
$M_1$ and $M_2$ are torsional submodules of $M$.
However, $M_1$ does not belong to $T(M_2)$ (note that $M_2$ is even a brick, thus it has no
non-trivial torsional submodules, see 4.3).
	\medskip
{\bf 9.7. Lemma.} 
{\it Let $B$ be a brick and $M$ a non-zero module in $\Cal E(B)$. Then
$M$ has an endomorphism with image a brick, 
and $B$ is the only brick which occurs in this way.} Thus, the type of a non-zero
homogeneous module is uniquely determined.
	\medskip
Proof. First, we show that $B$ occurs as the image of an endomorphism of $M$.
Since $M$ belongs to $\Cal E(B)$, there is a filtration
$0 = M_0 \subset M_1 \subset \cdots \subset M_m = M$ of $M$ such that all 
factors are isomorphic to $B$. A corresponding map $M \to M/M_{m-1} \simeq B \simeq M_1 
\subseteq M$ is an endomorphism of $M$ which image isomorphic to $B$. 

Conversely, let $f$ be an endomorphism of $M$ whose image is a brick. Since $\Cal E(B)$
is an exact abelian subcategory, the image $M'$ of $f$ belongs to $\Cal E(B)$.
Now $M'$ is a non-zero module in $\Cal E(B)$. As we have seen in the first part of the proof, $M'$ has an endomorphism with image $f(M')$ being isomorphic to $B$. 
But we assume that $M'$ is a brick, thus the image of an endomorphism of $M'$ is either
zero or $M'$ itself. This shows that $M' = f(M')$, thus $M'$ is isomorphic to $B$. 
$\s$ 
	\medskip
{\bf Remark.} Let $M$ be a homogeneous module. 
According to 9.7, the endomorphism
ring of $M$ shows whether $M$ is a brick or not.
But $\End(M)$ gives only limited information
about $M$. In particular, $\End(M)$ may be a $k$-algebra of dimension 2, whereas $M$
has a filtration with arbitrarily many factors of the form $B$. 
Here is an example. Consider the subring $A = k+J$ 
of the ring of $(t\times t)$-matrices with $t\ge 2,$
where $J$ is the set of nilpotent upper triangular matrices; 
and look at the set $M = k^t$ of column vectors.
The $A$-module $M$ is a serial module of length $t$, but
$\dim_k \End(M) = 2.$ 
	\bigskip
{\bf 10. Further remarks about brick chains and brick chain filtrations.}  
	\medskip

{\bf 10.1.} A filtration  $0 = M_0 \subseteq M_1 \subseteq \cdots \subseteq M_m = M$ 
will be said to be {\it solid} provided $\Hom(M_i/M_{i-1},M_j/M_{j-1}) = 0$ for
all $1 \le i < j \le m.$ 
	\medskip
{\bf Proposition.} {\it Let $M$ be a module with a filtration $0 = M_0 \subset M_1 
\subset \cdots \subset M_m = M.$ Then $(M_i)_i$ is a brick chain filtration iff 
$(M_i)_i$ is a solid filtration and all the factors are homogeneous.}
	\medskip
Proof. First, assume that $(M_i)_i$ is a brick chain filtration, say of type 
$(B_1,\dots,B_m).$ Since $M_i/M_{i-1}$ belongs to $\Cal E(B_i)$, all the factors of the
filtration are homogeneous. Also, for $i<j$, we have $\Hom(B_i,B_j) = 0.$ Therefore
$\Hom(M_i/M_{i-1},M_j/M_{j-1}) = 0$.

Conversely, assume that $(M_i)_i$ is a solid filtration (with proper
inclusions) and all factors are homogeneous. Now $F_i = M_i/M_{i-1}$ belongs to
$\Cal E(B_i)$ for some brick $B_i$. Since $F_i$
is non-zero, $B_i$ occurs both as a submodule and as a factor module of $F_i$.
Thus, any non-zero map $f\:B_i \to B_j$ yields a non-zero map 
$F_i \to F_j$. Since the given filtration is solid, we see that $\Hom(B_i,B_j) = 0$
for $i < j.$ Thus, $(B_1,\dots,B_m)$ is a brick chain. $\s$
	\medskip
{\bf 10.2.} The composition factors in the top of a module 
give rise to brick chain filtrations:
	\medskip
{\bf  Proposition.} {\it Let $M$ be a module. Let $S$ be a simple module which 
occurs in the top of $M$. Let $M'$ be the smallest submodule of $M$ such that $M/M'$
belongs to $\Cal E(S).$ Then, any torsional brick chain filtration $(M_i)_{0\le i\le m'}$ 
of $M'$ is part of a brick chain filtration $(M_i)_{0\le i\le m'+1}$ of $M$}
(with $M_{m'} = M'$).
	\smallskip
This shows: for any simple module $S$ in the top of $M$, there is 
a brick chain filtration whose type ends in $S$. Using duality, we also have:
for any simple module $S$ in the socle of $M$, there is 
a brick chain filtration whose type starts with $S$.
	\medskip
Proof. Let $(B_1,\dots,B_{m'})$ be the type of the torsional brick chain filtration of $M'.$
The bricks $B_i$ 
are in $\Cal T(M')$, thus the top of $B_i$ is generated by $M'.$ 
As a consequence, $\Hom(B_i,S) = 0$ for $1\le i \le m'.$ Thus
$(B_1,\dots,B_{m'},S)$ is a brick
chain, and $(M_i)_{1\le i \le {m'+1}}$ is a brick chain filtration
of type $(B_1,\dots,B_{m'},S).$
$\s$
	\medskip
{\bf 10.3. Proposition.} {\it A module $M$ has only one brick chain filtration iff
all composition factors of $M$ are isomorphic}. 
	\medskip
Proof. If all composition factors of $M$ are isomorphic, then clearly $(0 \subseteq M)$ is the
only brick chain filtration. 
Conversely, assume that $M$ has only one brick chain filtration.
According to 9.2, $M$ has only one top brick, say $B$.
On the other hand, 10.2 shows that $M$ has a brick chain filtration of 
type $(B_1,\dots,B_{m})$ with $B_{m} = S$ simple. This shows that $\Cal B(M) = \{S\}$.
Now 5.5 asserts that $T(M) = T(\Cal B(M))$. Since $M$ belongs to $T(M) = T(\Cal B(M)) = 
T(S)$, all comporisiton factors of $M$ are isomorphic to $S$.
$\s$
	\medskip
{\bf 10.4. Induced brick chain filtrations.} 
We say that a filtration $(M_i)_i$ of a module $M$ is {\it proper}
provided all the inclusions $M_{i-1} \subseteq M_i$ are proper. 
Of course, any filtration $(M_i)_i$
yields a proper filtration by deleting all the submodules $M_i$ with $M_{i-1} = M_i.$
Until  now, all the filtrations considered in the paper were proper. 
Let us call an arbitrary filtration of
a module $M$ a {\it brick chain filtration with repetitions} provided the corresponding
proper filtration is a brick chain filtration.
	\medskip
{\bf  Proposition.} (a) {\it Let $M,N$ be modules. Let $(X_i)_i$ be a 
brick chain filtration of  $X = M\oplus N$. Then 
$(X_i\cap M)_i$ is a brick chain filtration of $M$ with repetitions,} we say that it is
{\it induced from}
$(X_i)_i$.
	
(b) {\it Let $M$ be a module with a brick chain filtration 
$(M_i)_i$ say of type $(B_1,\dots,B_t)$. Let $N = \bigoplus B_i$.
Then, also the module $M\oplus N$ has a brick chain filtration of type 
$(B_1,\dots,B_t)$, this filtration is torsional, and it induces 
the given filtration $(M_i)_i$.}
This shows that {\it any brick chain filtration is induced from a torsional brick chain filtration.}
	\medskip
Proof. (a) Since $X_i$ is the torsion submodule of $X$ for the torsion class 
$T(B_1,\dots, B_i)$, we have $(X_i\cap M)\oplus (X_i\cap N) = X_i$.
Thus $(X_i\cap M)/(X_{i-1}\cap M)\oplus (X_i\cap N)/(X_{i-1}\cap N) =
X_i/X_{i-1}$ belongs to $\Cal E(B_i)$, thus $(X_i\cap M)/(X_{i-1}\cap M)$
belongs to $\Cal E(B_i)$. Of course, $(X_i\cap M)/(X_{i-1}\cap M)$ may be zero
(thus, the filtration may have repetitions).
	
(b) Let $N = \bigoplus_{j=1}^{t-1}B_j$ and $X = M\oplus N = X_t.$ 
For $0 \le i < t,$ let
$X_i = M_i\oplus \bigoplus_{j=1}^{i}B_j.$  
Then $X_t/X_{t-1} = M_t/M_{t-1} \in \Cal E(B_t).$ Also, 
for $1\le i < t,$ we have $X_i/X_{i-1} = M_i/M_{i-1}\oplus B_i\in \Cal E(B_i)$.
Thus, we see that $(X_i)_i$ is a brick chain filtration of type $(B_1,\dots,B_t).$

All the bricks $B_i$ with $1\le i < t$ are factor modules of $X,$
thus $T(B_1,\dots,B_{t-1}) \subseteq T(X)$. Since $X_{t-1}$ belongs to
$\Cal E(B_1,\dots,B_{t-1})$, we see that $X_{t-1}$ belongs to $T(X).$
Similarly, for $1\le i < t$, we have $T(B_1,\dots,B_{i-1}) \subseteq T(X_i).$
Again, $X_{i-1}$ belongs to $\Cal E(B_1,\dots,B_{i-1})$, thus
$X_{i-1}$ belongs to $T(X_i).$ This shows that $(X_i)_i$ is a torsional filtration.
$\s$
	\medskip
{\bf 10.5. Remark.} According to 10.4, brick chain filtrations of direct
sums yield brick chain filtrations of the summands. The converse is not true.
For example, let $A$ be a cyclic Nakayama algebra with two simple modules
and $M, N$ the two indecomposable modules of length two. Then $M, N$ are bricks, thus
they have brick chain filtrations of length one. But any brick chain filtration
of the module $M\oplus N$ has length three.
	\bigskip
\vfill\eject
{\bf 11. Appendix. Complete brick chains.}
	\medskip
{\bf 11.1. Definition.}
Following Demonet [De], one should also 
deal with arbitrary (not necessarily finite) brick chains: 
these are 
arbitrarily large totally ordered sets $\Cal B = \{B_i\mid i\in I\}$ of bricks with the
$\Hom$-condition $\Hom(B_i,B_j) = 0$ for all $i < j$.
Given brick chains $\Cal B, \Cal B'$, one calls $\Cal B$ a {\it refinement} of $\Cal B'$, provided
$\Cal B'$ is obtained from $\Cal B$ by deleting some elements. And, a brick chain is
said to be {\it complete} provided it has no proper refinement.
Similarly, dealing with arbitrary chains of torsion classes, there is the corresponding concept of refinements and of completeness of such chains.
	\medskip
Demonet has shown the following assertions.
	\smallskip
{\bf (1)} {\it Any brick chain has a complete refinement.}
	\smallskip
{\bf (2)} {\it The complete chains of torsion classes
correspond bijectively to the complete brick chains, sending a complete chain $(\Cal T_i)_i$ 
of torsion classes to the brick chain given by the labels of the
neighbor torsion classes in the chain $(\Cal T_i)_i$.}
	\smallskip
{\bf (3)} {\it Any torsion class is generated by the bricks of a (usually infinite) brick chain.}
	\smallskip
{\bf (4)} {\it Any simple module occurs in any complete brick chain.}
	\medskip
The proofs require to deal with sets of arbitrary cardinality, but otherwise they 
are easy. For (2), Corollary 7.4 is essential. Assertion (3) is a direct consequence of (2).
$\s$
	\bigskip
{\bf 11.2. Theorem.} {\it Let $\Cal B$ be a complete brick chain. 
Any module $M$ has a unique brick chain filtration 
$(M_i)_i$ such that $\Cal B$ is a refinement of the type of $(M_i)_i.$}
	\medskip
Proof, by induction on the length of $M$. The assertion is clear for the zero module. Thus, assume
that $M \neq 0.$ Let $(\Cal T_i)_{i\in I}$ be a complete chain of torsion classes.
Let $\Cal T = \bigcap_{M\in \Cal T_i} \Cal T_i$ and
$\Cal T' = \bigcup_{M\notin \Cal T_i} \Cal T_i$. 
Since we deal with a complete chain of torsion classes, both $\Cal T, \Cal T'$ belong to the chain.
Of course, we have $\Cal T'
\subseteq \Cal T$,  but even $\Cal T' \subset \Cal T$, since $M$ belongs to 
$\Cal T$ and does not belong to $\Cal T'$. 
It follows, that $\Cal T' \subset \Cal T$ are neighbors,
say with label $B\in \Cal B$, thus $\Cal T' = \Cal T\cap {}^\perp B.$ 

First, we show that $M$ has at most one brick chain such that the bricks of the type  
belong to $\Cal B$. 
Let $(M_i)_i$ be a brick chain of $M$, say of type $(B_1,\dots,B_t)$ and assume that all 
the bricks $B_i$ belong to $\Cal B.$ We want to show that $B_t = B.$
Not all the modules $B_i$ can belong to $\Cal T'$, since otherwise $M$ is in $\Cal T'$.
Thus, $B_t$ is not in $\Cal T'$. Now $M$ belongs to $\Cal T$ and $B_t$ is a factor module of $M$,
thus $B_t$ belongs to $\Cal T$. But $B$ is the only brick in $\Cal B$ which belongs to
$\Cal T \setminus \Cal T'$, therefore $B_t = B.$ By induction, also the bricks $B_1,\dots,B_{t-1}$
are determined by $M_{t-1}$ and $\Cal B$.

Conversely, we show that $M$ has a brick chain filtration say with type $\Cal B'$, 
such that $\Cal B$ is a refinement of $\Cal B'.$
Since $\Cal T = T(\Cal T',B) = T(\Cal T\cap{}^\perp B,B),$ 
we can apply 6.1 and obtain a submodule $M'$ of $M$ belonging
to $\Cal T'$ such that $M/M'$ is a non-zero module in $\Cal E(B)$. By induction, $M'$ has
a brick chain filtration say of type $\Cal B''$, where $\Cal B$ is a refinement of
$\Cal B''$.
$\s$
	\bigskip
{\bf 12. History and relevance (as well as additions).}
	\medskip
{\bf 12.1.} The results presented here are usually considered as part of the so-called
$\tau$-tilting theory (see 12.13). 
There is a strange reluctance
to deal with bricks. For example, many papers prefer to speak about 
$\tau$-tilting finiteness instead of brick finiteness, but these properties are
equivalent (see [DIJ]; here, $\tau$-tilting finiteness means that there
are only finitely many basic support-$\tau$-tilting modules:
In my opinion, brick finiteness 
is very easy to grasp, whereas $\tau$-tilting finiteness is much 
less intuitive). For our report, there was no need 
to mention $\tau$-tilting 
notions, nor even the Auslander-Reiten translation $\tau$ itself, thus we 
have avoided it.
In this way, we stress the completely elementary nature of the corresponding results.
Some remarks on $\tau$-tilting theory will be given in 12.13.

It is astonishing that the relevance of bricks when dealing 
with tilting modules and torsion classes, 
was observed only so late!
	\smallskip
{\bf 12.2. Bricks and semibricks.} The terminology ``semibrick'' is
due to Asai [A1].
I used to call a semibrick an ``antichain'' of bricks, 
but this is in conflict with Demonet's
important notion of a brick chain (to 
say that ``an antichain of bricks is a brick chain'', would sound rather odd).  
	\smallskip
{\bf 12.3. Torsion pairs $(\Cal T,\Cal F)$.} Torsion pairs 
were introduced by Dickson [Di] as a generalization of the use of torsion
and $p$-torsion subgroups of abelian groups, for dealing with  
arbitrary $R$-modules, were $R$ is any ring.
	\smallskip
{\bf 12.4. Hereditary torsion pairs, torsional submodules.}
In contrast to the classical example, torsion classes in general are not hereditary
(where {\it hereditary} means that the torsion class $\Cal T$ is closed under
submodules). The torsion classes $T(M)$ considered in our paper
are usually not hereditary.
But it turns out that dealing with a module $M$,
it is important to look at submodules
of $M$ which do belong to $T(M)$ (the torsional submodules). Our
focus on torsional submodules is 
an attempt to stress heredity properties for non-hereditary torsion classes. 

Theorems 1.2 and 3.2 should be seen in the light of the 
original example of abelian group theory: any finitely generated abelian group $M$
has a filtration $(M_i)_{0\le i \le m}$ 
where the first factors $M_i/M_{i-1}$ are in $\Cal E(\Bbb Z/p_i\Bbb Z),$
for pairwise different prime numbers $p_i$,  
whereas $M_m/M_{m-1}$ is in $\Cal E(\Bbb Z)$, and this filtration always splits! 
In our case, we cannot expect that
the filtrations in 1.2 and 3.2 
split. (It comes as a surprise that actually in first
examples one looks at, say dealing with Kronecker modules, many brick chain filtrations
do split.)
	\smallskip
{\bf 12.5. Auslander and Smal\o{} (and Demonet).} 
The relevance of torsion classes when dealing with finite length
categories was seen already by Auslander and Smal\o{} [AS]. 
		
Looking at a module category, the existence of cyclic paths in
the category or even in the Auslander-Reiten quiver, provides a lot of
difficulties. Only the 
representation-directed algebras are easy to visualize, but 
representation-directedness is a very special property. 
There have been many attempts to overcome the difficulties which
arise from the presence of cyclic paths. There is 
the covering theory by Gabriel and his school;
also, the book of  Auslander, Reiten, Smal\o{} is full of helpful devices:
to avoid short chains, to avoid short cycles.
However, all these methods are designed just for special, well-behaved situations.
If one wants to deal with an arbitrary module category,
the use of torsion classes always works. 
The reference to torsion classes allows to
consider the set of semibricks as
a partially ordered set.
 In this way, 
Demonet's proposal to look at brick chains stresses a very interesting
directedness feature of an arbitrary module category.
	\smallskip
{\bf 12.6. Wide subcategories and torsion classes.} Given an abelian category, 
the exact abelian subcategories which are closed under extensions
are now often called {\it wide} subcategories. 
The rather obvious relationship between semibricks and wide
subcategories was mentioned in [R1] under the name ``simplification''.
The search for semibricks (or wide subcategories) which
generate a given torsion class was initiated by Ingalls and Thomas [IT]. 
Theorem 2.3 generalizes some of their considerations.

The relevance of the endotop of a module is well-known
and was stressed by Asai when looking at $\tau$-rigid modules
(our proof of 5.5 follows closely Asai [A]). For a general study of 
widely generated torsion classes, see Asai and Pfeifer [AP] and 
Marks and Stovi\v cek [MS].
	\smallskip
{\bf 12.7. Homogeneous subcategories.}
The homogeneous subcategories are equivalent to the module category of a local 
algebra (not necessarily an artin algebra).
Actually, not much is known about the representation theory of 
local algebras which are not commutative! The commutative rings are 
studied very well in commutative algebra, but who cares about 
the non-commutative ones? Note that they can behave rather
differently and really deserve attention. 

For example, if $A$ is a commutative local ring, and $M$ 
is a serial module, then there is an endomorphism of $M$ with image
$\rad M$, thus $\et M$ is just the simple module. 
On the other hand, consider the subring $A = k+J$ of the ring of $(t\times t)$-matrices
where $J$ is the set of nilpotent upper triangular matrices,
as mentioned already in 5.2 (2).  This is a rather nice 
local ring; it is non-commutative provided $t \ge 3$.
The set $M = k^t$ of column vectors is a serial $A$-module. 
Here, the image of any non-invertible endomorphism of $M$ has length
at most one, thus we see that $\et M$ has dimension $t\!-\!1$;
in particular, it is not a brick provided $t\ge 3$
(for $t\ge 3$, we have $\et^\infty M = \et^{t-2} M = k \neq \et M$).
	\smallskip
{\bf 12.8. Neighbors of torsion classes.} 
Neighbor torsion classes $\Cal T' \subset \Cal T''$ 
have attracted a lot of interest and several
different denominations are used in the literature: 
that $\Cal T''$ covers $\Cal T'$, that there is
an arrow $\Cal T'' \to \Cal T'$ in the Hasse quiver of the lattice of torsion classes, 
or one speaks about minimal inclusions of torsion classes. 

As we have seen, it is easy to determine the lower neighbors of a finitely generated torsion class (and there are only finitely many),   
but unfortunately, it is difficult to deal with the upper neighbors: usually, there
may be infinitely many. 
For any torsion class $\Cal T$, the best way to find its upper neighbors seems to be
to look at the corresponding torsion free class $\Cal F$ and to try to determine
its lower neighbors, since the lower neighbors of $\Cal F$ correspond to the
upper neighbors of $\Cal T$.
	\smallskip
{\bf 12.9. Brick labeling.} The brick labeling as presented in section 7
was started for functorially finite
torsion classes in [AIR] and Asai [A] identified the labels as bricks.
The general case is due to Barnard, Carroll and Zhu [BCZ].
The brick $B$ used as label for the neighbor torsion classes $\Cal T' \subset \Cal T$ 
is called a {\it minimal extending module} for $\Cal T'$ in [BCZ].
In [AHL], the labels are said to be {\it torsion, 
nearly torsionfree} for the torsion pair $(\Cal T,\Cal T^\perp)$.
The bijection 2.9 between bricks and completely join irreducible torsion classes 
has been exhibited in Theorem 1.0.5 in [BCZ].
	\smallskip
{\bf 12.10. Brick chains.} Given chains of torsion classes, 
the brick labeling of the neighbor torsion classes yields the $\Hom$-condition 
which defines the brick chains. 
This observation was used by Demonet [De] to deal with
arbitrarily large totally ordered sets of bricks with this
$\Hom$-condition, see 10.8.
But note that already for 
the Kronecker algebra, the sets which occur explode:  
The Kronecker algebra $A$ over the field with 2 elements has cardinality 16, 
thus it is very easy to envision, but 
there are uncountably many 
complete brick chains (one is finite, all others are countable). 

Demonet's bijection 10.8 (2) is very charming: whereas torsion classes are complicated categories
(usually they are quite difficult to exhibit),
a brick is 
just a brick, a well-behaved single module! Of course, there is the next level: to look at
brick chains: they are not so easy to grasp, since usually these chains may be very large sets. But for 
getting a feeling for large brick chains, one may first restrict
the attention to finite brick chains (as we do in this report). Actually, in some cases,  there is no problem to overview well the complete brick: 
for example, an infinite complete brick chain for the Kronecker algebra 
has a preprojective and a preinjective part (both are easy to visualize);
and in-between, there is the regular part: it is 
$\Bbb P^1(k)$ with an arbitrary
(and really irrelevant) total ordering. 
	\smallskip
{\bf 12.11. Special brick chain filtrations} have been used already in [R2]: In this paper
we have shown that for a hereditary $k$-algebra, with
$k$ an algebraically closed field, any brick without self-extensions
is a tree module.
The basis of the proof is Schofield induction: it deals with certain brick chain filtrations
of length two. 
The brick chain filtrations used have type $(B_1,B_2)$,
where both $B_1,B_2$ are again bricks without self-extensions. 
These brick chain filtrations are never
torsional (since they are filtrations of bricks). 
This shows the
relevance of brick chain filtrations which are not torsional.
	\smallskip
{\bf 12.12. Artinian rings.} In this report, we have assumed to be in the context of artin algebras. 
Actually, nearly all the results
presented here are valid more generally in arbitrary length categories, thus
for finitely generated modules over left artinian rings. 
	 \smallskip
{\bf 12.13. Functorially finite torsion classes.}
Functorially finite torsion classes were first considered by Auslander 
and Smal\o{}. In 1984, Smal\o{} formulated the tie between functorially finite torsion classes
in $\mod A$  and tilting modules for suitable factor algebras of $A$. 
The basic objects of $\tau$-tilting theory are the $\tau$-rigid modules
(a module $M$ is {\it $\tau$-rigid} provided $\Hom(M,\tau M) = 0,$ where $\tau$
is the Auslander-Reiten translation): The 
Adachi-Iyama-Reiten paper [AIR] (published in 2014) showed that the functorially finite 
torsion classes are just the torsion classes generated by $\tau$-rigid modules,
a very important observation!
(Actually, several of the main results of [AIR] also follow from investigations by 
Derksen and Fei which were available in arXiv already in 2009, but it seems that this
paper is not really appreciated by the tau-tilting community).

We cannot give here even a concise summary of
$\tau$-tilting theory, but there do exist already many surveys which can be consulted.
Only few details should be mentioned: For an algebra with $n$ simple
modules, any functorially finite torsion class has precisely $n$, and
sufficiently many, neighbors, and all neighbors are again functorially finite.
Asai [A2] calls a torsion class {\it bicompact} provided
it has finitely many, and sufficiently many, neighbors, and he conjectures 
that a bicompact torsion class has to be functorially finite
(as we have mentioned, the converse is true). 
An important feature of $\tau$-tilting theory is the use of mutations in order
to describe for a given functorially finite torsion class its neighbors.
In the setting of tilting theory, mutations were studied by Happel and Unger (and others), 
but it took a long time 
that the relevance for arbitrary module categories was realized. 
	
We have seen in 2.5 that an algebra is brick finite iff it is torsion
class finite. It turns out that in this case, all torsion classes are not only finitely generated,
but even functorially finite. And conversely, if any torsion class is
functorially finite, then the algebra has to be brick finite. 
One can use this fact in order to show: if all torsion classes are 
finitely generated, then the algebra is brick finite, as mentioned in 2.5.
	\smallskip
{\bf 12.14. Further developments.} 
We say that a module $M$ has {\it 
brick chain complexity} at most $t$ provided there is a brick chain filtration of $M$ with $t$ 
factors. The {\it brick chain complexity} of an algebra $A$ is the 
supremum of the brick chain complexity of the indecomposable $A$-modules
(it is a natural number or $\infty$); it will be discussed in [R3[.
	\smallskip
{\bf Acknowledgment.} First of all,
the referee has to be praised for a careful reading of the
manuscript, finding a lot of awful inaccuracies and providing several very helpful
comments.
The author has to thank first Zhi-Wei Li 
and then the referee
for pointing out that one of the essential assertions, namely
Proposition 6.1, had to be reformulated twice. 
In addition, he is grateful to S\. Asai and W\. Crawley-Boevey
for useful remarks. 
	\bigskip
{\bf 13. References.}
	\smallskip
\item{[AHL]} L\. Angeleri H\"ugel, 
   I\. Herzog, R\. Laking. Simples in a cotilting heart. 
    Math\. Z\. 307, 12 (2024). 
\item{[AIR]} T\. Adachi, O\. Iyama, I\. Reiten. $\tau$-tilting theory. 
    Compos. Math., 150(3), 415--452 (2014). 
\item{[A1]} S\. Asai. Semibricks. Int. Math. Res. Not., (16), 4993--5054.
    \newline
    doi.org/10.1093/imrn/rny150 (2020).
\item{[A2]} S\. Asai. Bicompact torsion classes and conjectures on brick infinite algebras.
  \newline   arXiv:2604.04505
\item{[AP]} S\. Asai, C\. Pfeifer. 
   Wide subcategories and lattices of torsion classes.
   Algebras and Representation Theory.
   Vol. 25, 1611--1629 (2022). arXiv:1905.01148
\item{[AS]} M\. Auslander, S\. O\. Smal\o. 
    Addendum to: Almost split sequences in subcategories. J.
    Algebra, 71, 592--594 (1981).
\item{[BCZ]} E\. Barnard, A\. Carroll, S\. Zhu. Minimal inclusions of 
            torsion classes. Algebr. Comb., 2(5), 879--901 (2019).
\item{[De]} L\. Demonet. Maximal chains of torsion classes. Appendix to:
     B\. Keller. A survey on maximal green sequences,
     In: Representation theory and beyond, 
    Contemp. Math. 758, Amer. Math. Soc., Providence, RI, 267–-286
     (2020). 
\item{[Di]} S\. E\. Dickson. A torsion theory for Abelian categories. 
     Trans. Amer. Math. Soc., 121, 223--235 (1966).
\item{[DIJ]} L\. Demonet, O\. Iyama, G\. Jasso. 
    $\tau$--tilting finite algebras, bricks, and g-vectors. Int.
    Math. Res. Not. IMRN, (3), 852--892 (2019).
\item{[IT]} C\. Ingalls, H\. Thomas. 
    Noncrossing partitions and representations of quivers.
    Comp. Math. 145, 1533–156 (2009).
\item{[MS]} F\. Marks, J\. Stovi\v cek. Torsion classes, wide subcategories and
    localisations. Bulletin of the London Mathematical Society.
    Bull. Lond. Math. Soc. 49, 405--416, doi:10.1112/blms.12033  (2017).
\item{[R1]} C\. M\. Ringel. Representations of $K$-species and bimodules.
     J. Algebra 41, 269--302 (1976).
\item{[R2]} C\. M\. Ringel. Exceptional modules are tree modules.
     Lin. Alg. Appl. 275-276, 471--493  (1998). 
\item{[R3]} C\. M\. Ringel. 
     The brick chain complexity of an artin algebra. Preprint. 

    	\medskip
{\baselineskip=1pt
\rmk
Claus Michael Ringel\par
Fakult\"at f\"ur Mathematik, Universit\"at Bielefeld \par
POBox 100131, D-33501 Bielefeld, Germany  \par
ringel\@math.uni-bielefeld.de}

\bye